\documentclass[11pt]{article}

\usepackage[a4paper,margin=1in]{geometry}
\usepackage{graphicx}
\usepackage{booktabs}
\usepackage{multirow}
\usepackage{xcolor}
\usepackage{textcomp}

\usepackage{amsmath,amssymb,amsfonts}
\usepackage{amsthm}
\usepackage{mathrsfs}

\usepackage{dsfont}
\newcommand{\Law}{\mathcal{L}}
\usepackage{bbm}

\usepackage{algorithm}
\usepackage{algorithmicx}
\usepackage{algpseudocode}
\usepackage{listings}

\usepackage[hidelinks]{hyperref}
\usepackage[nameinlink,capitalise]{cleveref}

\graphicspath{{./}{./figures/}}

\newcommand{\rev}[1]{#1}
\newcommand{\revtitle}[1]{#1}

\newtheorem{theorem}{Theorem}
\newtheorem{proposition}[theorem]{Proposition}
\newtheorem{lemma}[theorem]{Lemma}
\newtheorem{corollary}[theorem]{Corollary}

\theoremstyle{definition}
\newtheorem{definition}{Definition}

\theoremstyle{remark}
\newtheorem{remark}{Remark}

\title{On Lorden's Inequality and Renewal-Type Processes with Dependent Inter-renewal Times}

\author{
El'mira Yu. Kalimulina\thanks{Research Department, Incarnet Math Modelling, 0012, 32 Komitas Ave, Yerevan, Armenia. Moscow State University (MSU), Leninskie Gory, Moscow, 119991, Russia.  Institute for Information Transmission Problems of RAS, Bolshoy Karetny per., Moscow, 127051, Russia. \texttt{elmira.yu.k@gmail.com}} \and
Galina A. Zverkina\thanks{Institute for Information Transmission Problems of RAS, Bolshoy Karetny per., Moscow, 127051, Russia.}
}
\date{} 

\begin{document}
\maketitle
\begin{abstract}
\rev{We consider renewal-type processes whose positive inter-renewal times may be dependent, non-identically distributed, and may have mixed distributions. We introduce a generalised intensity measure extending the classical hazard-rate representation to this setting. Under a two-sided comparison scheme for the inter-renewal laws and an additional renewal-measure domination condition \textnormal{(RD)}, we prove a Lorden-type bound for the forward recurrence time. This bound provides an explicit first-moment input for coupling constructions and, once the remaining coupling parameters are controlled, yields a total-variation estimate. We illustrate the result on exponential, mixed, Markov-modulated, and Pareto benchmarks. In the i.i.d.\ benchmarks, the bound has the correct renewal scale up to a universal factor; in the Markov-modulated benchmark, the explicit Lorden constant is verified while the final convergence consequence remains conditional on \textnormal{(RD)}; and in the Pareto case the construction identifies the natural second-moment threshold for finiteness of the Lorden input.}
\end{abstract}


\section{Introduction}

\rev{A basic problem in queueing theory and reliability theory is to obtain quantitative bounds on convergence to stationarity, rather than only existence and uniqueness of an invariant distribution. In many models the embedded renewal structure is non-classical: inter-renewal times may be dependent, need not be identically distributed, and may have mixed laws. In such settings, generator-based, spectral, or Lyapunov methods are often model-specific and need not yield explicit computable bounds.}

\rev{We study renewal-type processes with positive inter-renewal times under a comparison scheme based on lower and upper laws for the inter-renewal variables. The main theorem, Theorem~\ref{thm:generalized_lorden}, proves a Lorden-type bound for the forward recurrence time under Assumptions~\ref{ass:grp:min}--\ref{ass:grp:moment} together with the additional renewal-measure domination condition~\eqref{cond:RD}. The bound provides an explicit first-moment input for the coupling constructions recorded in Sections~\ref{sec:stateOftheArt} and~\ref{section:classicalresults}. Corollary~\ref{cor:generalized_lorden_coupling} then gives only a conditional coupling consequence: once a stationary version, a trial coupling scheme, and the remaining coupling parameters are available, one obtains a total-variation estimate.}

\rev{The benchmark calculations serve different purposes. The exponential and mixed i.i.d.\ examples calibrate the constant against the classical renewal scale. The Markov-modulated example shows that the comparison assumptions and the explicit constant $C_{\mathrm{GL}}^{\mathrm{MME}}$ can be verified in closed form under dependence and heterogeneity, but it does not verify the full setwise domination condition~\eqref{cond:RD}; the corresponding convergence-rate statement therefore remains conditional. The Pareto example identifies the natural threshold $\alpha>2$ for finiteness of the first Lorden input.}

\rev{The paper is organised as follows. Section~\ref{sec:stateOftheArt} recalls the classical renewal and coupling material used later. Section~\ref{section:classicalresults} records the age-process coupling scheme in the classical setting. Section~\ref{sec:mainresultLordenGener} states and proves the generalised Lorden-type inequality and its conditional coupling consequence. Section~\ref{sec:application} contains the benchmark examples, while Appendix~\ref{app:benchmark_tables} collects the comparative tables.}

\section{Classical Lorden's Inequality and Methods for Estimating the Convergence Rate for Regenerative Processes}\label{sec:stateOftheArt}

\rev{Lorden's inequality is a classical uniform bound for the excess over a boundary in a renewal process with nonnegative i.i.d.\ inter-renewal times \cite{Lorden1970,CarlssonNerman1986,Chang1994}. Applications to convergence-rate questions for regenerative models are discussed, for instance, in \cite{zverkina2020lordensinequalitypolynomialrate}. In this section we fix the renewal-theoretic terminology used later and recall the classical Lorden and coupling material needed in the sequel; see, e.g., \cite{1958_Smith,cox1962renewal,Feller71,asmussen2003applied,zverkina2020lordensinequalitypolynomialrate}.}

\subsection{Lorden's Inequality for the Classical Renewal Process}\label{subsec:stateOftheArtLordenClassic}

We begin with the standard framework. Throughout this subsection we work in the ordinary (non-delayed) setting.

\begin{definition}[Regenerative process, cf.\ \cite{asmussen2003applied,Borovkov2020}]
A stochastic process \(X=(X_t)_{t\geq 0}\) is called \emph{regenerative} if there exists an increasing sequence of almost surely finite random times
\[
0=T_0<T_1<T_2<\cdots,\qquad T_n\to\infty \quad \text{a.s.},
\]
such that the cycle segments
\[
\left(\, \bigl(X_{T_{n-1}+s}\bigr)_{0\leq s<T_n-T_{n-1}},\; T_n-T_{n-1}\right),\qquad n\geq 1,
\]
are independent and identically distributed. The times \(T_n\) are called \emph{regeneration times}.
\end{definition}
\begin{definition}[Ordinary renewal process]
Let \(\xi_1,\xi_2,\ldots\) be i.i.d.\ nonnegative random variables with distribution function
\[
F(x)=\mathbb{P}(\xi_1\leq x),\qquad x\geq 0.
\]
Define
\[
S_0:=0,\qquad S_n:=\sum_{j=1}^n \xi_j,\qquad n\geq 1.
\]
The counting process
\[
N_t:=\sup\{n\geq 0:S_n\leq t\},\qquad t\geq 0,
\]
is called the \emph{renewal process} generated by \(\{\xi_n\}_{n\geq 1}\), and the points \(S_n\) are its \emph{renewal epochs}.
\end{definition}

\begin{definition}[Embedded renewal process]
Let \(X\) be a regenerative process with regeneration times \(\{T_n\}_{n\geq 0}\). Then the cycle lengths
\[
\xi_n:=T_n-T_{n-1},\qquad n\geq 1,
\]
form an i.i.d.\ sequence of nonnegative random variables, and the corresponding renewal process is called the \emph{embedded renewal process} associated with \(X\). In this case the renewal epochs coincide with the regeneration times:
\[
S_n=T_n,\qquad n\geq 0.
\]
When convenient, we also write \(t_n:=S_n\) for the renewal epochs.
\end{definition}

\begin{definition}[Backward recurrence time]
Let \(N=(N_t)_{t\geq 0}\) be a renewal process with renewal epochs \(\{S_n\}_{n\geq 0}\). The \emph{backward recurrence time} (or \emph{age}) at time \(t\geq 0\) is
\[
B_t:=t-S_{N_t}.
\]
Equivalently, if one writes \(t_n:=S_n\), then \(B_t=t-t_{N_t}\).
\end{definition}

\begin{definition}[Forward recurrence time]
Let \(N=(N_t)_{t\geq 0}\) be a renewal process with renewal epochs \(\{S_n\}_{n\geq 0}\). The \emph{forward recurrence time} (or \emph{excess}, or \emph{residual life}) at time \(t\geq 0\) is
\[
W_t:=S_{N_t+1}-t.
\]
Equivalently, if one writes \(t_n:=S_n\), then \(W_t=t_{N_t+1}-t\).
\end{definition}

Throughout the paper we use the convention
\[
N_t=\sup\{n\geq 0:S_n\leq t\}.
\]
With this convention, if
\[
\nu(b):=\inf\{n\geq 1:S_n>b\},\qquad b\geq 0,
\]
then the overshoot over the level \(b\) is
\[
R_b:=S_{\nu(b)}-b,
\]
and, for deterministic \(t\geq 0\),
\[
R_t=W_t.
\]
Accordingly, the classical Lorden inequality is a statement about the \emph{forward recurrence time} \(W_t\) (equivalently, the overshoot \(R_t\)), not about the age \(B_t\). The corresponding undershoot at level \(b\) is \(b-S_{\nu(b)-1}\), which in the renewal notation equals \(B_b\).

\begin{theorem}[Lorden's inequality \cite{Lorden1970}]\label{teorem:Lorden}
Let \(\xi_1,\xi_2,\ldots\) be i.i.d.\ nonnegative random variables such that
\[
0<\mu:=\mathbb{E}\xi_1<\infty,
\qquad
\mathbb{E}\xi_1^{\,2}<\infty.
\]
Define
\[
S_0:=0,\qquad S_n:=\sum_{j=1}^n \xi_j,\qquad
\nu(b):=\inf\{n\geq 1:S_n>b\},
\qquad
R_b:=S_{\nu(b)}-b.
\]
Then, for every \(b\geq 0\),
\[
\mathbb{E}R_b\leq \frac{\mathbb{E}\xi_1^{\,2}}{\mathbb{E}\xi_1}.
\]
Equivalently, since \(R_t=W_t\) under the convention \(N_t=\sup\{n\geq 0:S_n\leq t\}\),
\[
\sup_{t\geq 0}\mathbb{E}W_t\leq \frac{\mathbb{E}\xi_1^{\,2}}{\mathbb{E}\xi_1}.
\]
\end{theorem}

\begin{remark}
Alternative proofs of Theorem~\ref{teorem:Lorden} are given in \cite{CarlssonNerman1986}. Higher-moment analogues were obtained by Chang \cite{Chang1994}; see also \cite{Sugakova2007}. Extensions beyond the classical i.i.d.\ setting include independent non-identically distributed summands \cite{Spouge2007} and random environments \cite{Svensson02}. For a recent discussion of overshoot asymptotics beyond the nonnegative renewal setting, see \cite{Wong2025}. In the present subsection we use only the classical first-moment form stated in Theorem~\ref{teorem:Lorden}.
\end{remark}

Assume that \(\xi_1,\xi_2,\dots\) are i.i.d., \(\xi_1>0\) almost surely,
\[
0<\mu:=\mathbb{E}\xi_1<\infty,
\]
and that the distribution \(F\) is non-arithmetic. Then classical renewal theory implies that both the age \(B_t\) and the excess \(W_t\) converge in distribution to the common equilibrium distribution
\[
H_e(x):=\frac{1}{\mu}\int_0^x \overline F(u)\,du,
\qquad x\ge 0,
\]
where \(\overline F(x)=1-F(x)\); see \cite{1958_Smith,cox1962renewal,Feller71,asmussen2003applied}. More precisely,
\[
\mathbb{P}(B_t\le x)\longrightarrow H_e(x),
\qquad
\mathbb{P}(W_t\le x)\longrightarrow H_e(x),
\qquad x\ge 0,
\quad t\to\infty.
\]
Equivalently,
\[
\mathbb{P}(W_t>x)\longrightarrow \frac{1}{\mu}\int_x^\infty \overline F(u)\,du,
\qquad x\ge 0.
\]
Moreover, if \(\mathbb{E}\xi_1^{\,2}<\infty\), then the renewal equations for \(b(t):=\mathbb{E}B_t\) and \(w(t):=\mathbb{E}W_t\), together with the key renewal theorem, imply that
\[
\mathbb{E}B_t\longrightarrow \frac{\mathbb{E}\xi_1^{\,2}}{2\mu},
\qquad
\mathbb{E}W_t\longrightarrow \frac{\mathbb{E}\xi_1^{\,2}}{2\mu},
\qquad t\to\infty;
\]
see, for example, \cite{1958_Smith,cox1962renewal,Feller71,asmussen2003applied}.
These asymptotic relations describe only the equilibrium behaviour
of \(B_t\) and \(W_t\) as \(t\to\infty\).
They do not imply the existence of a uniform finite-time bound of
Lorden type. Such bounds require additional arguments.

\subsection{Convergence Rate Estimates for Regenerative Processes}\label{subsec:convergenceIntroPart}

In the regenerative setting, convergence rates are obtained by coupling a process started from a given initial state with a stationary version. The rate is then controlled by the tail of the coupling time. The basic tool used below is the following standard coupling estimate; see, for example, \cite{lindvall2002lectures,Thorisson00}.

\begin{proposition}[Standard coupling estimate]\label{prop:coupling_stationary}
Let \(X=(X_t)_{t\ge 0}\) be a stochastic process on a measurable state space \((E,\mathcal E)\), and let \(\pi\) be a stationary distribution. For each initial state \(x\in E\), suppose that one can construct on a common probability space
\begin{itemize}
    \item a version \(X^x=(X_t^x)_{t\ge 0}\) of the process started from \(x\),
    \item a stationary version \(\widetilde X=(\widetilde X_t)_{t\ge 0}\) of the same process, whose one-time marginal law is \(\pi\),
    \item and a nonnegative almost surely finite random variable \(\tau_x\),
\end{itemize}
such that
\[
X_t^x=\widetilde X_t \qquad \text{for all } t\ge \tau_x \quad \text{almost surely}.
\]
Then, for every \(t\ge 0\),
\[
\bigl\|\Law(X_t^x)-\pi\bigr\|_{\mathrm{TV}}
\le \mathbb{P}(\tau_x>t).
\]
More generally, if \(\varphi:[0,\infty)\to[0,\infty)\) is measurable, non-decreasing, and \(\varphi(t)>0\) for \(t>0\), and if
\[
\mathbb{E}\,\varphi(\tau_x)<\infty,
\]
then
\[
\bigl\|\Law(X_t^x)-\pi\bigr\|_{\mathrm{TV}}
\le
\frac{\mathbb{E}\,\varphi(\tau_x)}{\varphi(t)},
\qquad t>0.
\]
In particular:
\begin{enumerate}
    \item if \(\mathbb{E}\tau_x^\alpha<\infty\) for some \(\alpha>0\), then for every \(0<\beta\le \alpha\),
    \[
    \bigl\|\Law(X_t^x)-\pi\bigr\|_{\mathrm{TV}}
    \le
    \frac{K_\beta(x)}{t^\beta},
    \qquad
    K_\beta(x):=\mathbb{E}\tau_x^\beta;
    \]
    \item if \(\mathbb{E}e^{a\tau_x}<\infty\) for some \(a>0\), then for every \(0<\beta\le a\),
    \[
    \bigl\|\Law(X_t^x)-\pi\bigr\|_{\mathrm{TV}}
    \le
    K_{\exp,\beta}(x)e^{-\beta t},
    \qquad
    K_{\exp,\beta}(x):=\mathbb{E}e^{\beta\tau_x}.
    \]
\end{enumerate}
\end{proposition}

Proposition~\ref{prop:coupling_stationary} is purely abstract: it reduces convergence-rate estimates to control of a coalescence time, but does not by itself guarantee existence of a stationary regime. In the regenerative setting, existence and convergence require separate assumptions on the regeneration structure.

More generally, if \(X\) and \(X'\) are defined on a common probability space and satisfy
\[
X_t=X_t' \qquad \text{for all } t\ge \tau \quad \text{almost surely},
\]
for some almost surely finite random time \(\tau\), then the same argument gives, for every measurable set \(S\subseteq E\) and every \(t>0\),
\[
\left|\mathbb P\{X_t\in S\}-\mathbb P\{X_t'\in S\}\right|
\le \mathbb P\{\tau>t\}
\le \frac{\mathbb E[\psi(\tau)]}{\psi(t)},
\]
whenever \(\psi:[0,\infty)\to[0,\infty)\) is measurable, non-decreasing, satisfies \(\psi(t)>0\) for \(t>0\), and \(\mathbb E[\psi(\tau)]<\infty\).

For regenerative processes with i.i.d.\ cycle lengths, the remaining problem is therefore to estimate the relevant coalescence time through renewal quantities. In the constructions considered later, Lorden-type bounds will be used precisely for this purpose.

\subsubsection{Coupling Method and a Common-Part Coupling Lemma}

For the concrete constructions used later, we shall need a standard common-part coupling lemma.

\begin{lemma}[Fundamental coupling lemma, cf.\ \cite{denHollander2012}]\label{BLC}
Let \(\xi_1\) and \(\xi_2\) have densities \(\phi_1\) and \(\phi_2\) with respect to Lebesgue measure on \(\mathbb R\), and assume that
\[
\varkappa:=\int_{-\infty}^{\infty}\min\{\phi_1(s),\phi_2(s)\}\,ds>0.
\]
Then there exist random variables \(\vartheta_1\) and \(\vartheta_2\) on a suitable probability space such that
\[
\vartheta_i\overset{D}{=}\xi_i,\qquad i=1,2,
\]
and
\[
\mathbb P\{\vartheta_1=\vartheta_2\} = \varkappa.
\]
\end{lemma}

\begin{proof}
Let
\[
\phi(s):=\min\{\phi_1(s),\phi_2(s)\},\qquad \varkappa=\int_{\mathbb R}\phi(s)\,ds.
\]
If \(\varkappa=1\), then \(\phi_1=\phi_2\) almost everywhere, and one may take \(\vartheta_1=\vartheta_2\) with that common distribution.

Assume now that \(0<\varkappa<1\). Define
\[
g(s):=\frac{\phi(s)}{\varkappa},
\qquad
h_i(s):=\frac{\phi_i(s)-\phi(s)}{1-\varkappa},
\qquad i=1,2.
\]
Then \(g,h_1,h_2\) are probability densities on \(\mathbb R\). Let \(B\) be a Bernoulli random variable with \(\mathbb P(B=1)=\varkappa\), let \(Y\) have density \(g\), and let \(Z_1,Z_2\) have densities \(h_1,h_2\), respectively, where \(B,Y,Z_1,Z_2\) are independent. Define
\[
\vartheta_i :=
\begin{cases}
Y, & \text{if } B=1,\\
Z_i, & \text{if } B=0,
\end{cases}
\qquad i=1,2.
\]
Then \(\vartheta_i\) has density
\[
\varkappa g +(1-\varkappa)h_i = \phi_i,
\]
so \(\vartheta_i\overset{D}{=}\xi_i\). Moreover, on \(\{B=1\}\) one has \(\vartheta_1=\vartheta_2=Y\), while on \(\{B=0\}\) the variables \(Z_1\) and \(Z_2\) are independent and absolutely continuous, hence \(\mathbb P\{Z_1=Z_2\}=0\). Therefore
\[
\mathbb P\{\vartheta_1=\vartheta_2\}=\varkappa.
\]
\end{proof}

\begin{lemma}[Generalisation of the coupling lemma]\label{lem:basecouplinglemma:n}
Let \(f_i\) be probability densities of random variables \(\theta_i\), \(i=1,\ldots,n\), and assume that
\[
\varkappa:=\int_{-\infty}^{\infty}\min_{1\le i\le n} f_i(s)\,ds>0.
\]
Then one can construct random variables \(\vartheta_1,\ldots,\vartheta_n\) on a suitable probability space such that
\[
\vartheta_i\overset{D}{=}\theta_i,\qquad i=1,\ldots,n,
\]
and
\[
\mathbb P\{\vartheta_1=\vartheta_2=\cdots=\vartheta_n\}=\varkappa.
\]
\end{lemma}

\begin{proof}
Let
\[
f(s):=\min_{1\le i\le n} f_i(s).
\]
If \(\varkappa=1\), then all \(f_i\) coincide almost everywhere, and one may take \(\vartheta_1=\cdots=\vartheta_n\) with this common distribution. Assume \(0<\varkappa<1\), and define
\[
g(s):=\frac{f(s)}{\varkappa},
\qquad
h_i(s):=\frac{f_i(s)-f(s)}{1-\varkappa},
\qquad i=1,\ldots,n.
\]
Then \(g,h_1,\ldots,h_n\) are probability densities. Let \(B\) be Bernoulli with \(\mathbb P(B=1)=\varkappa\), let \(Y\) have density \(g\), and let \(Z_1,\ldots,Z_n\) be independent random variables with densities \(h_1,\ldots,h_n\), independent of \(B\) and \(Y\). Set
\[
\vartheta_i :=
\begin{cases}
Y, & \text{if } B=1,\\
Z_i, & \text{if } B=0,
\end{cases}
\qquad i=1,\ldots,n.
\]
Then each \(\vartheta_i\) has density \(f_i\), and on \(\{B=1\}\) all variables are equal. On \(\{B=0\}\), the random variables \(Z_1,\ldots,Z_n\) are independent and absolutely continuous, hence
\[
\mathbb P\{Z_1=\cdots=Z_n\}=0.
\]
Therefore
\[
\mathbb P\{\vartheta_1=\cdots=\vartheta_n\}=\varkappa.
\]
\end{proof}

These lemmas provide the basic common-part coupling step used later in the construction of successful couplings for regenerative and renewal-type processes.

\subsubsection{Successful coupling}\label{sec:parallcoupl}

Let \(X^{(1)}=(X_t^{(1)})_{t\ge0}\) and \(X^{(2)}=(X_t^{(2)})_{t\ge0}\) be two versions of the same process started from \(x_1,x_2\in E\). A pair
\[
\mathcal Z=(Z_t^{(1)},Z_t^{(2)})_{t\ge0}
\]
is called a \emph{successful coupling} of \(X^{(1)}\) and \(X^{(2)}\) if
\[
\Law\bigl((Z_t^{(i)})_{t\ge0}\bigr)=\Law\bigl((X_t^{(i)})_{t\ge0}\bigr),
\qquad i=1,2,
\]
and there exists a nonnegative almost surely finite random variable \(\tau(x_1,x_2)\) such that
\[
Z_t^{(1)}=Z_t^{(2)}
\qquad \text{for all } t\ge \tau(x_1,x_2)
\quad \text{almost surely}.
\]

Therefore, by Proposition~\ref{prop:coupling_stationary},
\[
\|\Law(X_t^{(1)})-\Law(X_t^{(2)})\|_{\mathrm{TV}}
\le
\mathbb P\{\tau(x_1,x_2)>t\}
\le
\frac{\mathbb E\varphi(\tau(x_1,x_2))}{\varphi(t)},
\qquad t>0,
\]
whenever \(\varphi\) is non-decreasing, \(\varphi(t)>0\) for \(t>0\), and \(\mathbb E\varphi(\tau(x_1,x_2))<\infty\).

Hence, with
\[
\mathscr R(x_1,x_2):=\mathbb E\varphi(\tau(x_1,x_2)),
\]
one has
\[
\|\Law(X_t^{(1)})-\Law(X_t^{(2)})\|_{\mathrm{TV}}
\le
\frac{\mathscr R(x_1,x_2)}{\varphi(t)},
\qquad t>0.
\]

\section{Coupling framework for regenerative processes via the backward recurrence time}\label{section:classicalresults}

This section records the classical prototype underlying the later Lorden-based coupling argument. In particular, we isolate the transfer step from the age process to the regenerative process, describe the local transition mechanism of the classical age process, and formulate a schematic trial-coupling estimate used later.

\subsection{Distribution of the regenerative process via the backward recurrence time}

Let \(X=(X_t)_{t\ge0}\) be a regenerative process with state space \(E\), and let
\[
B_t:=t-S_{N_t}
\]
be the backward recurrence time associated with i.i.d.\ inter-renewal times
\(\xi_1,\xi_2,\ldots\) with distribution function \(F\).

\begin{lemma}[Equilibrium limit for the backward recurrence time]
\label{lem:weakconvergenceovershoot}
Assume that \(\xi_1>0\) almost surely, \(\mu:=\mathbb E\xi_1<\infty\), and the renewal law is non-arithmetic. Then
\[
B_t \Rightarrow \widetilde B
\qquad \text{as } t\to\infty,
\]
where
\[
\mathbb P\{\widetilde B\le s\}
=
\frac1\mu\int_0^s (1-F(u))\,du,
\qquad s\ge0.
\]
\end{lemma}

\begin{proof}
This is the classical equilibrium-age limit; see Subsection~\ref{subsec:stateOftheArtLordenClassic}.
\end{proof}

\begin{lemma}[Transfer of total-variation bounds from \(B_t\) to \(X_t\)]
\label{lem:sxodim}
Assume that there exists a time-independent Markov kernel
\begin{equation}\label{eq:ConditionK}
K:[0,\infty)\times\mathcal B(E)\to[0,1]
\end{equation}
such that, for every \(t\ge0\),
\[
\mathcal L(X_t)=\mathcal L(B_t)K.
\]
Let \(\widetilde B\) be as in Lemma~\ref{lem:weakconvergenceovershoot}, and define
\[
\pi_X:=\mathcal L(\widetilde B)K.
\]
If
\[
\|\mathcal L(B_t)-\mathcal L(\widetilde B)\|_{\mathrm{TV}}
\le d(t),
\qquad t\ge0,
\]
for some function \(d(t)\to0\), then
\[
\|\mathcal L(X_t)-\pi_X\|_{\mathrm{TV}}
\le d(t),
\qquad t\ge0.
\]
In particular,
\[
\mathcal L(X_t)\Rightarrow \pi_X.
\]
\end{lemma}

\begin{proof}
Write \(\mu_t=\mathcal L(B_t)\) and \(\widetilde\mu=\mathcal L(\widetilde B)\).
Then \(\mathcal L(X_t)=\mu_tK\) and \(\pi_X=\widetilde\mu K\).
Since Markov kernels are contractions in total variation,
\[
\|\mu_tK-\widetilde\mu K\|_{\mathrm{TV}}
\le
\|\mu_t-\widetilde\mu\|_{\mathrm{TV}},
\]
which yields the claim.
\end{proof}

\begin{remark}
Condition~\eqref{eq:ConditionK} holds, for example, if \(X_t=h(B_t)\) for a measurable function \(h\), or more generally if the conditional law of \(X_t\) given \(B_t\) does not depend on \(t\). In applications, one must verify that \(\pi_X\) is indeed the marginal law of a stationary version of \(X\).
\end{remark}

\subsection{Hazard rate and transition mechanism of the age process}

We restrict attention here to the classical absolutely continuous case. Let
\[
\xi_1,\xi_2,\ldots
\]
be i.i.d.\ positive inter-renewal times with distribution function \(F\), density \(f\), and survival function
\[
\overline F(s):=1-F(s),
\qquad s\ge0.
\]
Set
\[
E_B:=\{y\ge0:\ \overline F(y)>0\},
\]
and, writing \(t_n:=S_n\) for the renewal epochs,
\[
B_t:=t-t_{N_t}.
\]

For \(y\in E_B\), the residual time to the next renewal conditional on \(B_t=y\) has distribution
\[
\mathbb P\{W_t>u\mid B_t=y\}
=
\frac{\overline F(y+u)}{\overline F(y)},
\qquad u\ge0.
\]
Hence, for every \(\Delta>0\),
\[
\mathbb P\{\text{at least one renewal on }(t,t+\Delta]\mid B_t=y\}
=
1-\frac{\overline F(y+\Delta)}{\overline F(y)}.
\]

Define the hazard rate by
\begin{equation}\label{eq:hazardRateDefinition1}
\lambda(y):=\frac{f(y)}{\overline F(y)},
\qquad y\in E_B.
\end{equation}
Then
\begin{equation}\label{eq:hazardRateDefinition2}
\frac{\overline F(y+\Delta)}{\overline F(y)}
=
\exp\!\left(-\int_y^{y+\Delta}\lambda(u)\,du\right),
\end{equation}
and therefore
\begin{equation}\label{eq:hazardRateDefinition3}
\mathbb P\{\text{at least one renewal on }(t,t+\Delta]\mid B_t=y\}
=
1-\exp\!\left(-\int_y^{y+\Delta}\lambda(u)\,du\right).
\end{equation}
If, in addition, \(\lambda\) is continuous at \(y\), then
\begin{equation}\label{eq:hazardRateDefinition4FolkDefinition}
\mathbb P\{\text{at least one renewal on }(t,t+\Delta]\mid B_t=y\}
=
\lambda(y)\Delta+o(\Delta),
\qquad \Delta\downarrow0.
\end{equation}

Equivalently,
\[
\overline F(s)=\exp\!\left(-\int_0^s\lambda(u)\,du\right),
\qquad
f(s)=\lambda(s)\exp\!\left(-\int_0^s\lambda(u)\,du\right).
\]

The age process \(B_t\) is Markov. More precisely, for \(y\in E_B\),
\[
P_\Delta(y,A)
=
\frac{\overline F(y+\Delta)}{\overline F(y)}\,\mathbf 1_{\{y+\Delta\in A\}}
+
\int_0^\Delta
\mathbb P_0\{B_{\Delta-r}\in A\}\,
\frac{f(y+r)}{\overline F(y)}\,dr,
\qquad
A\in\mathcal B([0,\infty)),
\]
where \(P_\Delta(y,\cdot)\) denotes the conditional law of \(B_{t+\Delta}\) given \(B_t=y\), and \(\mathbb P_0\) is the law of the age process started from \(0\).

In particular, if \(g\in C_b^1([0,\infty))\) and \(\lambda\) is locally bounded and continuous, then
\[
\mathbb E[g(B_{t+\Delta})\mid B_t=y]
=
g(y)+\Delta\bigl(g'(y)+\lambda(y)(g(0)-g(y))\bigr)+o(\Delta),
\qquad \Delta\downarrow0.
\]
Thus the local operator is
\[
\mathcal A g(y):=g'(y)+\lambda(y)\bigl(g(0)-g(y)\bigr),
\qquad y\in E_B,
\]
in the sense of the above first-order expansion.

For mixed inter-renewal laws, the density-based hazard description is no longer sufficient and must be replaced by a transition-kernel or hazard-measure formulation; this is done in Section~\ref{sec:mainresultLordenGener}.

\subsection{Convergence rate estimates for the age process in the classical renewal setting}

For \(b\in E_B\), let \(R_0^{\,b}\) be a nonnegative random variable such that
\[
\mathbb P\{R_0^{\,b}>u\}
=
\frac{\overline F(b+u)}{\overline F(b)},
\qquad u\ge0,
\]
and assume that \(R_0^{\,b}\) is independent of an ordinary age process
\(\bar B=(\bar B_t)_{t\ge0}\) started from \(0\) and driven by the inter-renewal law \(F\).
Define
\[
B_0^b:=b,
\qquad
B_t^b:=
\begin{cases}
b+t, & 0<t<R_0^{\,b},\\[1mm]
\bar B_{\,t-R_0^{\,b}}, & t\ge R_0^{\,b},
\end{cases}
\qquad
\mathcal P_t^b:=\Law(B_t^b).
\]

Assume \(\mu=\mathbb E\xi_1<\infty\), and let \(\pi_B\) be the equilibrium law of Lemma~\ref{lem:weakconvergenceovershoot}. When convergence from the ordinary initial condition is invoked, we additionally assume that the renewal law is non-arithmetic.

\begin{theorem}[Coupling-based rate estimate for the classical age process]
\label{thm:classical_B_rate}
Suppose that for some \(b\in E_B\) one can construct on a common probability space a coupling of \(B^b=(B_t^b)_{t\ge0}\) with a stationary version \(\widetilde B=(\widetilde B_t)_{t\ge0}\), whose one-time marginal law is \(\pi_B\), together with a coalescence time \(\tau_b\) such that
\[
B_t^b=\widetilde B_t
\qquad \text{for all } t\ge\tau_b
\quad \text{almost surely}.
\]
Then, for every \(t\ge0\),
\[
\|\mathcal P_t^b-\pi_B\|_{\mathrm{TV}}
\le
\mathbb P\{\tau_b>t\}.
\]
More generally, if \(\varphi:[0,\infty)\to[0,\infty)\) is measurable, non-decreasing, \(\varphi(t)>0\) for \(t>0\), and
\[
\mathbb E\varphi(\tau_b)<\infty,
\]
then
\[
\|\mathcal P_t^b-\pi_B\|_{\mathrm{TV}}
\le
\frac{\mathbb E\varphi(\tau_b)}{\varphi(t)},
\qquad t>0.
\]
In particular:
\begin{enumerate}
    \item if \(\mathbb E\tau_b^{\,k-1}<\infty\) for some \(k\ge2\), then
    \[
    \|\mathcal P_t^b-\pi_B\|_{\mathrm{TV}}
    \le
    \frac{\mathbb E\tau_b^{\,k-1}}{t^{k-1}},
    \qquad t>0;
    \]
    \item if \(\mathbb E e^{\beta\tau_b}<\infty\) for some \(\beta>0\), then
    \[
    \|\mathcal P_t^b-\pi_B\|_{\mathrm{TV}}
    \le
    \mathbb E e^{\beta\tau_b}\,e^{-\beta t},
    \qquad t>0.
    \]
\end{enumerate}
\end{theorem}

\begin{proof}
This is Proposition~\ref{prop:coupling_stationary} applied to \(B^b\).
\end{proof}

Thus the convergence-rate problem reduces to constructing a coupling and estimating its coalescence time.

\subsubsection{Parallel coupling construction for the classical age process}

In the classical absolutely continuous setting, the coupling step is formulated at the level of residual-life distributions. If the current age is \(u\in E_B\), then the residual time to the next renewal has density
\[
f_u(s):=\frac{f(u+s)}{\overline F(u)},
\qquad s\ge0.
\]
Accordingly, common-part couplings are governed by overlap coefficients
\[
\int_0^\infty \min\{f_u(s),f_v(s)\}\,ds,
\]
as in Lemma~\ref{BLC}. A detailed construction in the classical i.i.d.\ case is given in \cite{ZverkinaDCCN2017}. Here we isolate only the abstract trial-coupling template used later.

\begin{proposition}[A schematic trial-coupling estimate for the classical age process]
\label{prop:trial_coupling_classical}
Fix \(b\in E_B\). Let \(B^b=(B_t^b)_{t\ge0}\) be the classical age process started from \(b\), and let \(\widetilde B=(\widetilde B_t)_{t\ge0}\) be a stationary version with one-time marginal law \(\pi_B\). Suppose that on a filtered probability space
\[
(\Omega,\mathcal F,(\mathcal F_t)_{t\ge0},\mathbb P)
\]
one can construct a coupling of \(B^b\) and \(\widetilde B\), together with an increasing sequence of almost surely finite stopping times
\[
0=T_0<T_1<T_2<\cdots
\]
and events \(A_n\in\mathcal F_{T_n}\), \(n\ge1\), such that:
\begin{enumerate}
    \item \(\mathbb E T_1\le a_b\) for some \(a_b<\infty\);
    \item for every \(n\ge1\),
    \[
    \mathbb E[T_{n+1}-T_n\mid \mathcal F_{T_n}] \le m
    \qquad \text{a.s.},
    \]
    for some \(m<\infty\);
    \item for every \(n\ge1\),
    \[
    \mathbb P(A_n\mid \mathcal F_{T_{n-1}})\ge q_b
    \qquad \text{a.s. on } \bigcap_{k=1}^{n-1}A_k^c,
    \]
    for some \(q_b\in(0,1]\);
    \item on \(A_n\), the coupled processes coalesce from time \(T_n\) onward:
    \[
    B_t^b=\widetilde B_t
    \qquad \text{for all } t\ge T_n
    \quad \text{a.s. on } A_n.
    \]
\end{enumerate}
Define
\[
N:=\inf\{n\ge1:\ A_n \text{ occurs}\}.
\]
Then \(N<\infty\) almost surely, \(T_N\) is a coalescence time, and
\[
\mathbb P(N>n)\le (1-q_b)^n,
\qquad n\ge0.
\]
Moreover,
\[
\mathbb E T_N\le a_b+m\left(\frac{1}{q_b}-1\right).
\]
Consequently,
\[
\|\mathcal P_t^b-\pi_B\|_{\mathrm{TV}}
\le
\frac{a_b+m(q_b^{-1}-1)}{t},
\qquad t>0.
\]
\end{proposition}

\begin{proof}
Set
\[
\mathcal G_n:=\mathcal F_{T_n},
\qquad n\ge0.
\]
Then
\[
\{N>n\}=\bigcap_{k=1}^n A_k^c\in\mathcal G_n,
\qquad n\ge0.
\]
On \(\{N=n\}=A_n\cap\bigcap_{k=1}^{n-1}A_k^c\), assumption~(iv) implies
\[
B_t^b=\widetilde B_t
\qquad \text{for all } t\ge T_n,
\]
so \(T_N\) is a coalescence time whenever \(N<\infty\).

For \(n\ge1\),
\[
\mathbb P(N>n)
=
\mathbb E\!\left[
\mathbf 1_{\{N>n-1\}}
\mathbb E(\mathbf 1_{A_n^c}\mid\mathcal G_{n-1})
\right].
\]
Since \(\{N>n-1\}\in\mathcal G_{n-1}\) and, on \(\{N>n-1\}\),
\[
\mathbb E(\mathbf 1_{A_n^c}\mid\mathcal G_{n-1})
=
1-\mathbb P(A_n\mid\mathcal G_{n-1})
\le 1-q_b,
\]
it follows that
\[
\mathbb P(N>n)\le(1-q_b)\mathbb P(N>n-1).
\]
Hence
\[
\mathbb P(N>n)\le(1-q_b)^n,
\qquad n\ge0,
\]
and therefore \(N<\infty\) almost surely.

Since \(T_n\) is increasing and \(N<\infty\) a.s.,
\[
T_N
=
T_1+\sum_{n\ge1}(T_{n+1}-T_n)\mathbf 1_{\{N>n\}}
\qquad \text{a.s.}
\]
By Tonelli's theorem,
\[
\mathbb E T_N
=
\mathbb E T_1
+
\sum_{n\ge1}\mathbb E\!\left[(T_{n+1}-T_n)\mathbf 1_{\{N>n\}}\right].
\]
Because \(\{N>n\}\in\mathcal G_n\), assumption~(ii) yields
\[
\mathbb E\!\left[(T_{n+1}-T_n)\mathbf 1_{\{N>n\}}\right]
=
\mathbb E\!\left[
\mathbf 1_{\{N>n\}}
\mathbb E(T_{n+1}-T_n\mid\mathcal G_n)
\right]
\le m\,\mathbb P(N>n).
\]
Therefore,
\[
\mathbb E T_N
\le
a_b+m\sum_{n\ge1}(1-q_b)^n
=
a_b+m\left(\frac{1}{q_b}-1\right).
\]
The final estimate follows from Proposition~\ref{prop:coupling_stationary} with \(\varphi(u)=u\).
\end{proof}

\begin{corollary}[Classical Lorden input for the trial-coupling template]
\label{cor:lorden_input}
In the setting of Proposition~\ref{prop:trial_coupling_classical}, assume in addition that \(\mathbb E\xi_1^{\,2}<\infty\), and suppose that the trial times \(T_n\) are chosen at successive renewal epochs of the stationary component \(\widetilde B\). Then
\[
a_b:=\mathbb E T_1
=
\mathbb E \widetilde W_0
=
\frac{\mathbb E\xi_1^{\,2}}{2\mu}
\le
\frac{\mathbb E\xi_1^{\,2}}{\mu},
\]
where \(\widetilde W_0\) is the stationary forward recurrence time at time \(0\). Thus the initial-delay term in Proposition~\ref{prop:trial_coupling_classical} is controlled by the classical Lorden constant.
\end{corollary}

\begin{proof}
Under the stated choice of trial times,
\[
T_1\overset{D}{=}\widetilde W_0.
\]
The equilibrium mean identity for the forward recurrence time yields
\[
\mathbb E T_1
=
\mathbb E \widetilde W_0
=
\frac{\mathbb E\xi_1^{\,2}}{2\mu},
\]
and the last inequality is immediate.
\end{proof}

\begin{remark}[How Lorden enters the classical coupling construction]
\label{rem:classical_coupling_mechanism}
Proposition~\ref{prop:trial_coupling_classical} separates two inputs. First, one needs a lower bound on the success probability of a single coupling attempt, expressed in the absolutely continuous case through overlap coefficients of residual-life densities:
\[
\varkappa_\Theta
:=
\inf_{0\le u\le \Theta}
\int_0^\infty \min\{f_u(s),f(s)\}\,ds.
\]
If \(\varkappa_\Theta>0\), then Lemma~\ref{BLC} yields a successful coupling whenever the current age lies in \([0,\Theta]\). Second, one needs control of the trial times. If trials are made at renewal epochs of the stationary component, then Corollary~\ref{cor:lorden_input} bounds the first trial delay by the classical Lorden constant.
\end{remark}

\begin{remark}[On higher-order rates]
The first-moment estimate above yields a \(t^{-1}\)-bound. Higher polynomial or exponential convergence rates require corresponding higher-moment or exponential bounds for the coalescence time.
\end{remark}

\section{Generalised Lorden inequality and coupling consequences}\label{sec:mainresultLordenGener}

Sections~\ref{sec:stateOftheArt}--\ref{section:classicalresults} developed the classical renewal-based coupling scheme. Its two key ingredients are a bound on the forward recurrence time, which controls the initial delay in the coupling construction, and a coupling step based on overlap of residual-life laws. In the absolutely continuous i.i.d.\ setting, these residual-life laws are expressed through the hazard rate.

The purpose of the present section is to extend that mechanism to a renewal-type setting in which the inter-renewal times may be dependent, non-identically distributed, and may contain both absolutely continuous and atomic parts. The main technical device is a generalised intensity measure, which replaces the density-based hazard while preserving the comparison structure needed for the later Lorden-type and coupling arguments.

Throughout this section, the notation for the counting process, the age process, and the forward recurrence time is consistent with Sections~\ref{subsec:stateOftheArtLordenClassic} and~\ref{section:classicalresults}.

\subsection{Renewal-type setting}

Let \(\{\xi_i\}_{i\ge1}\) be positive random variables defined on a common probability space. Set
\[
S_0:=0,
\qquad
S_n:=\sum_{j=1}^n \xi_j,
\qquad n\ge1,
\]
and assume that
\[
S_n\to\infty
\qquad \text{almost surely as } n\to\infty.
\]
Define the associated counting process by
\[
N_t:=\sup\{n\ge0:S_n\le t\},
\qquad t\ge0,
\]
and write
\[
F_i(x):=\mathbb P\{\xi_i\le x\},
\qquad x\ge0,\quad i\ge1.
\]

The variables \(\xi_i\) may be dependent and need not be identically distributed. The main objects used below are the generalised intensity measures of the inter-renewal laws and the comparison structure imposed by the standing assumptions introduced later in this section.

\subsubsection{Generalised intensity}

Since the inter-renewal laws considered below may contain both absolutely continuous and atomic parts, the density-based hazard used in Section~\ref{section:classicalresults} is no longer sufficient. We therefore work with a measure-valued hazard object, called the \emph{generalised intensity measure}; see, for example, \cite{borovkov2022compound}.

Throughout this subsection, a positive random variable will be called \emph{mixed} if its law has no singular continuous part, equivalently, if it consists only of an absolutely continuous part and an at most countable atomic part.

\begin{definition}[Generalised intensity measure]
Let \(\xi\) be a positive mixed random variable with distribution function \(F\), survival function
\[
\overline F(x):=1-F(x),
\qquad x\ge0,
\]
density \(f\) of the absolutely continuous part of \(F\), and set of atoms
\[
A_F:=\{a>0:\Delta F(a):=F(a)-F(a-)>0\}.
\]
The \emph{generalised intensity measure} of \(\xi\) is the extended nonnegative Borel measure \(\lambda_\xi\) on \([0,\infty)\) defined by
\[
\lambda_\xi(ds)
:=
\frac{f(s)}{\overline F(s-)}\,ds
+
\sum_{a\in A_F}
\left(
-\log\frac{\overline F(a)}{\overline F(a-)}
\right)\delta_a(ds),
\]
where \(f(s)/\overline F(s-)\) is understood as \(0\) whenever \(\overline F(s-)=0\), and the atomic coefficient is understood as \(+\infty\) if \(\overline F(a)=0\). Equivalently, if \(\overline F(a)=0\), then
\[
\lambda_\xi((0,x])=\infty
\qquad \text{for all } x\ge a.
\]
\end{definition}

If \(\overline F(s-)=0\), then \(f(s)=0\) for Lebesgue-a.e.\ such \(s\), since the absolutely continuous part of \(F\) has no mass beyond a point where \(F\) has already reached \(1\). Thus the convention \(0/0:=0\) is consistent.

In the purely absolutely continuous case,
\[
\lambda_\xi(ds)=\frac{f(s)}{\overline F(s)}\,ds,
\]
so the definition reduces to the usual hazard-rate representation.

\begin{proposition}[Reconstruction from the generalised intensity measure]
\label{prop:reconstruction_generalised_intensity}
For the measure \(\lambda_\xi\) defined above, the survival function of \(\xi\) is given by
\begin{equation}\label{eq:generalized_intensity_reconstruction}
\mathbb P\{\xi>x\}
=
\exp\!\bigl(-\lambda_\xi((0,x])\bigr),
\qquad x\ge0.
\end{equation}
Equivalently,
\[
F(x)=1-\exp\!\bigl(-\lambda_\xi((0,x])\bigr),
\qquad x\ge0.
\]
In particular, the distribution of \(\xi\) is uniquely determined by its generalised intensity measure.
\end{proposition}

\begin{proof}
Let
\[
S(x):=\mathbb P\{\xi>x\}=\overline F(x),
\qquad x\ge0,
\]
and define
\[
x_*:=\inf\{x\ge0:S(x)=0\}\in(0,\infty].
\]

Fix \(x<x_*\). Then \(S\) is strictly positive on \([0,x]\). Since \(\xi\) is mixed, the càdlàg nonincreasing function \(S\) has on \([0,x]\) only an absolutely continuous part and a jump part. Because \(u\mapsto \log u\) is \(C^1\) on \([S(x),1]\), the chain rule for functions of bounded variation yields
\[
d(\log S)(s)
=
-\frac{f(s)}{S(s-)}\,ds
+
\sum_{a\in A_F\cap(0,x]}
\log\frac{S(a)}{S(a-)}\,\delta_a(ds)
\]
on \((0,x]\). Integrating from \(0\) to \(x\) and using \(S(0)=1\) (since \(\xi>0\) almost surely), we obtain
\[
-\log S(x)
=
\int_{(0,x]} \frac{f(s)}{S(s-)}\,ds
+
\sum_{\substack{a\in A_F\\ a\le x}}
\left(-\log\frac{S(a)}{S(a-)}\right)
=
\lambda_\xi((0,x]).
\]
Exponentiating gives
\[
S(x)=\exp\!\bigl(-\lambda_\xi((0,x])\bigr),
\qquad x<x_*.
\]

If \(x\ge x_*\) and \(x_*<\infty\), then \(S(x)=0\). If \(S(x_*-) >0\) and \(S(x_*)=0\), then \(x_*\) is a terminal atom, so the corresponding atomic contribution in \(\lambda_\xi\) is \(+\infty\); hence
\[
\lambda_\xi((0,x])=\infty
\qquad \text{for all } x\ge x_*.
\]
If instead \(S(y)\downarrow0\) continuously as \(y\uparrow x_*\), then by the already established formula for \(y<x_*\),
\[
\lambda_\xi((0,y])=-\log S(y)\to\infty
\qquad \text{as } y\uparrow x_*,
\]
and again \(\lambda_\xi((0,x])=\infty\) for all \(x\ge x_*\). Thus \eqref{eq:generalized_intensity_reconstruction} holds for all \(x\ge0\), and the formula for \(F(x)\) follows immediately.
\end{proof}

\subsubsection{Assumptions on the distributions of recovery periods}\label{subsec:predp}

For the remainder of this section, we impose the following standing assumptions on the inter-renewal sequence \(\{\xi_i\}_{i\ge1}\).

\begin{enumerate}
    \item\label{ass:grp:min}
    For each \(i\ge1\),
    \[
    \xi_i=\min\{\zeta_i,\theta_i\},
    \]
    where \(\{\zeta_i\}_{i\ge1}\) are i.i.d.\ positive mixed random variables with common generalised intensity measure \(\varphi\), and \(\{\theta_i\}_{i\ge1}\) are positive mixed random variables with generalised intensity measures \(\mu_i\). The sequence \(\{\zeta_i\}_{i\ge1}\) is independent of the family \(\{\theta_i\}_{i\ge1}\).

    \item\label{ass:grp:upper}
    There exists an extended nonnegative Borel measure \(Q\) on \([0,\infty)\) such that
    \[
    \varphi+\mu_i\le Q
    \qquad \text{for all } i\ge1,
    \]
    in the sense of measure domination.

    \item\label{ass:grp:proper}
    The common generalised intensity measure \(\varphi\) satisfies
    \[
    \varphi((0,\infty))=\infty.
    \]

    \item\label{ass:grp:positivity}
    The absolutely continuous part of \(\varphi\) has a density, denoted by \(\varphi^{\mathrm{ac}}(s)\), such that
    \[
    \varphi^{\mathrm{ac}}(s)>0
    \qquad \text{for a.e. } s>T
    \]
    for some \(T\ge0\).

    \item\label{ass:grp:moment}
    For some \(k\ge2\),
    \[
    \int_0^\infty x^{k-1}\exp\!\bigl(-\varphi((0,x])\bigr)\,dx<\infty.
    \]
\end{enumerate}

No independence is assumed within the family \(\{\theta_i\}_{i\ge1}\), so the sequence \(\{\xi_i\}_{i\ge1}\) may be dependent. Assumption~\ref{ass:grp:proper} is automatic for the generalised intensity measure of a proper positive random variable, but is recorded explicitly because it will be used when defining the comparison law \(\Phi\). Since
\[
\varphi\le \varphi+\mu_i\le Q,
\]
Assumption~\ref{ass:grp:proper} also implies
\[
Q((0,\infty))=\infty,
\]
so the comparison law associated with \(Q\) is proper.

\begin{definition}[Quasi-renewal process]\label{def:quasiRenewal}
A renewal-type counting process satisfying Assumptions~\ref{ass:grp:min}--\ref{ass:grp:moment} will be called a \emph{quasi-renewal process}. This terminology is used only as a shorthand within the present paper.
\end{definition}

\begin{lemma}\label{lemma:min}
Under Assumption~\ref{ass:grp:min}, each \(\xi_i\) is a positive mixed random variable, and its generalised intensity measure satisfies
\[
\lambda_i=\varphi+\mu_i,
\qquad i\ge1.
\]
In particular, under Assumption~\ref{ass:grp:upper},
\[
\varphi\le \lambda_i\le Q,
\qquad i\ge1,
\]
in the sense of measure domination.
\end{lemma}

\begin{proof}
By Proposition~\ref{prop:reconstruction_generalised_intensity},
\[
\mathbb P\{\zeta_i>x\}
=
\exp\!\bigl(-\varphi((0,x])\bigr),
\qquad
\mathbb P\{\theta_i>x\}
=
\exp\!\bigl(-\mu_i((0,x])\bigr).
\]
Since \(\zeta_i\) and \(\theta_i\) are independent,
\[
\mathbb P\{\xi_i>x\}
=
\mathbb P\{\zeta_i>x,\ \theta_i>x\}
=
\exp\!\bigl(-(\varphi+\mu_i)((0,x])\bigr).
\]
Because \(\varphi+\mu_i\) again has only an absolutely continuous and an atomic part, Proposition~\ref{prop:reconstruction_generalised_intensity} implies that \(\xi_i\) is mixed and that its generalised intensity measure is \(\lambda_i=\varphi+\mu_i\). The inequalities
\[
\varphi\le\lambda_i\le Q
\]
follow immediately.
\end{proof}

Define the comparison distribution functions
\[
G(x):=1-\exp\!\bigl(-Q((0,x])\bigr),
\qquad
\Phi(x):=1-\exp\!\bigl(-\varphi((0,x])\bigr),
\qquad x\ge0.
\]
By Proposition~\ref{prop:reconstruction_generalised_intensity}, \(\Phi\) is the common distribution function of the i.i.d.\ variables \(\zeta_i\), and \(G\) is the proper comparison distribution associated with \(Q\).

These formulas are the measure-valued counterparts of the classical hazard-rate reconstruction formulas from Section~\ref{section:classicalresults}. In the purely absolutely continuous case, if
\[
Q(ds)=Q(s)\,ds,
\qquad
\varphi(ds)=\varphi(s)\,ds,
\]
then
\[
G(x)=1-\exp\!\left(-\int_0^x Q(s)\,ds\right),
\qquad
\Phi(x)=1-\exp\!\left(-\int_0^x \varphi(s)\,ds\right),
\]
which agrees with the classical formulas.

At this stage one obtains only one-dimensional comparison laws. In the dependent setting, such marginal bounds do not by themselves imply comparison of the partial sums \(S_n\) with sums of independent comparison variables. The lower convolution estimate below is available because Assumption~\ref{ass:grp:min} yields the pathwise inequality \(\xi_i\le\zeta_i\).

\begin{remark}\label{Remark:usl}
The following consequences of Assumptions~\ref{ass:grp:min}--\ref{ass:grp:moment} will be used later.
\begin{enumerate}
    \item For every \(i\ge1\) and every \(x\ge0\),
    \[
    G(x)\ge F_i(x)\ge \Phi(x).
    \]
    Equivalently, if \(\xi^{-}\) and \(\xi^{+}\) have distribution functions \(G\) and \(\Phi\), respectively, then
    \[
    \xi^{-}\le_{\mathrm{st}}\xi_i\le_{\mathrm{st}}\xi^{+},
    \qquad i\ge1.
    \]

    \item By Assumption~\ref{ass:grp:moment} and the tail-integral formula,
    \[
    \mathbb E(\xi^{+})^{k}
    =
    k\int_0^\infty x^{k-1}\exp\!\bigl(-\varphi((0,x])\bigr)\,dx<\infty.
    \]
    Hence
    \[
    \mathbb E\xi_i^{\,m}\le \mathbb E(\xi^{+})^{m}<\infty,
    \qquad 1\le m\le k,\quad i\ge1.
    \]

    \item Assumption~\ref{ass:grp:positivity} is a nondegeneracy condition on the common part \(\varphi\) used later in the coupling construction. Assumption~\ref{ass:grp:moment} is only a moment condition and does not by itself imply exponential convergence.
\end{enumerate}
\end{remark}

\begin{lemma}\label{lem:convolutioninequlatity}
Under Assumption~\ref{ass:grp:min}, for every \(n\ge1\),
\[
S_n\le \sum_{i=1}^n \zeta_i
\qquad \text{almost surely}.
\]
Consequently,
\[
\mathbb P\{S_n\le t\}
\ge
\Phi^{*n}(t),
\qquad t\ge0,
\]
where \(\Phi^{*n}\) denotes the \(n\)-fold convolution of \(\Phi\).
\end{lemma}

\begin{proof}
By Assumption~\ref{ass:grp:min},
\[
\xi_i=\min\{\zeta_i,\theta_i\}\le\zeta_i
\qquad \text{a.s. for every } i\ge1.
\]
Summing over \(i=1,\dots,n\) gives
\[
S_n\le \sum_{i=1}^n\zeta_i
\qquad \text{a.s.}
\]
Therefore
\[
\left\{\sum_{i=1}^n \zeta_i\le t\right\}\subseteq \{S_n\le t\},
\]
and hence
\[
\mathbb P\{S_n\le t\}\ge \Phi^{*n}(t).
\]
\end{proof}

\begin{remark}[On the scope of the convolution comparison]
The lower bound in Lemma~\ref{lem:convolutioninequlatity} is of pathwise origin. By contrast, the marginal inequalities
\[
F_i(x)\le G(x),\qquad x\ge0,
\]
do not imply an upper bound of the form
\[
\mathbb P\{S_n\le t\}\le G^{*n}(t)
\]
under arbitrary dependence of \(\xi_1,\dots,\xi_n\), since \(G^{*n}\) is the distribution function of a sum of independent comparison variables. Any upper comparison used later must therefore rely on additional structure.
\end{remark}

\subsection{Generalised Lorden's inequality}

As in Sections~\ref{subsec:stateOftheArtLordenClassic} and~\ref{section:classicalresults}, the object of interest is the forward recurrence time
\[
W_t:=S_{N_t+1}-t,
\qquad t\ge0.
\]
The structural assumptions of Subsection~\ref{subsec:predp} yield one-dimensional comparison laws and the lower convolution estimate of Lemma~\ref{lem:convolutioninequlatity}, but they do not provide an upper comparison for the partial sums \(S_n\). To obtain a uniform upper bound for \(W_t\), we therefore impose an additional hypothesis at the level of the associated renewal measures.

Let \(\zeta\) denote a generic copy of the i.i.d.\ comparison variable with distribution function \(\Phi\), and let \(\eta\) be a positive random variable with distribution function \(G\). By Remark~\ref{Remark:usl},
\[
\eta\le_{\mathrm{st}}\xi_i\le_{\mathrm{st}}\zeta,
\qquad i\ge1.
\]
Moreover,
\[
0<\mathbb E\eta\le \mathbb E\zeta<\infty,
\qquad
\mathbb E\eta^2\le \mathbb E\zeta^2<\infty.
\]
Set
\[
m_G:=\mathbb E\eta,
\qquad
L_G:=\frac{\mathbb E\eta^2}{(\mathbb E\eta)^2}.
\]

Let \(\eta_1,\eta_2,\ldots\) be i.i.d.\ copies of \(\eta\), and define
\[
T_0:=0,
\qquad
T_n:=\sum_{j=1}^n \eta_j,
\qquad n\ge1.
\]
The corresponding classical renewal measure is
\[
U_G(A):=\sum_{n=0}^\infty \mathbb P\{T_n\in A\},
\qquad A\in\mathcal B([0,\infty)).
\]
Likewise, for the renewal-type partial sums \(S_n\),
\[
U_\xi(A):=\sum_{n=0}^\infty \mathbb P\{S_n\in A\},
\qquad A\in\mathcal B([0,\infty)).
\]

\medskip
\noindent
\textbf{Condition (RD).}
We assume that
\begin{equation}\tag{RD}\label{cond:RD}
U_\xi(A)\le U_G(A),
\qquad
A\in\mathcal B([0,\infty)).
\end{equation}
This is an additional structural assumption. It does not follow from the marginal bounds \(F_i\le G\) alone and must be verified separately in applications.

\begin{lemma}[Classical renewal convolution bound]\label{lem:renewal_measure_convolution}
Let \(g:[0,\infty)\to[0,\infty)\) be nonnegative, nonincreasing, right-continuous, and integrable. Then, for every \(t\ge0\),
\[
\int_{[0,t]} g(t-s)\,U_G(ds)
\le
L_G\,g(0)+\frac{1}{m_G}\int_0^\infty g(u)\,du.
\]
\end{lemma}

\begin{proof}
Let
\[
N_t^G:=\sup\{n\ge0:T_n\le t\},
\qquad
W_t^G:=T_{N_t^G+1}-t,
\qquad t\ge0,
\]
be the counting process and forward recurrence time of the classical renewal process with inter-renewal law \(G\).

By Theorem~\ref{teorem:Lorden},
\[
\sup_{t\ge0}\mathbb E W_t^G
\le
\frac{\mathbb E\eta^2}{\mathbb E\eta}
=
L_G\,m_G.
\]
Hence the standard renewal identity gives
\[
U_G([0,t])
=
\mathbb E(N_t^G+1)
=
\frac{t+\mathbb E W_t^G}{m_G}
\le
\frac{t}{m_G}+L_G,
\qquad t\ge0.
\]

We next derive the increment bound
\begin{equation}\label{eq:renewal_increment_bound}
U_G((a,b]\cap[0,\infty))
\le
L_G+\frac{b-a}{m_G},
\qquad a\le b.
\end{equation}
If \(a<0\), then
\[
U_G((a,b]\cap[0,\infty))=U_G([0,b])
\le
L_G+\frac{b}{m_G}
\le
L_G+\frac{b-a}{m_G}.
\]
Assume now that \(0\le a\le b\), and write \(h:=b-a\). Conditionally on \(W_a^G\),
\[
N_b^G-N_a^G=
\mathbf 1_{\{W_a^G\le h\}}
\bigl(1+\widetilde N_{h-W_a^G}^G\bigr),
\]
where \(\widetilde N^G\) is an independent copy of \(N^G\). Therefore
\[
\mathbb E[N_b^G-N_a^G\mid W_a^G]
=
\mathbf 1_{\{W_a^G\le h\}}\,U_G([0,h-W_a^G])
\le
\sup_{0\le u\le h}U_G([0,u])
\le
L_G+\frac{h}{m_G}.
\]
Taking expectations yields \eqref{eq:renewal_increment_bound}.

For \(y>0\), define
\[
r(y):=\sup\{u\ge0:\ g(u)>y\}\in[0,\infty),
\]
with the convention \(r(y)=0\) if the set is empty. Since \(g\) is nonincreasing and right-continuous, the super-level set \(\{u\ge0:\ g(u)>y\}\) is the interval \([0,r(y))\), and \(r(y)<\infty\) for every \(y>0\).

Using the layer-cake representation and Tonelli's theorem,
\begin{align*}
\int_{[0,t]} g(t-s)\,U_G(ds)
&=
\int_{[0,t]}\int_0^\infty \mathbf 1_{\{g(t-s)>y\}}\,dy\,U_G(ds) \\
&=
\int_0^\infty U_G\bigl((t-r(y),t]\cap[0,\infty)\bigr)\,dy \\
&\le
\int_0^{g(0)}\left(L_G+\frac{r(y)}{m_G}\right)\,dy \\
&=
L_G\,g(0)+\frac{1}{m_G}\int_0^\infty r(y)\,dy.
\end{align*}
A second layer-cake identity gives
\[
\int_0^\infty r(y)\,dy=\int_0^\infty g(u)\,du,
\]
which proves the claim.
\end{proof}

\rev{The next theorem is the main result. It gives a Lorden-type upper bound for the forward recurrence time under Assumptions~\ref{ass:grp:min}--\ref{ass:grp:moment} and Condition~\eqref{cond:RD}.}

\begin{theorem}[Generalised Lorden-type inequality under renewal-measure domination]
\label{thm:generalized_lorden}
Assume Assumptions~\ref{ass:grp:min}--\ref{ass:grp:moment} and Condition~\eqref{cond:RD}. Then, for every \(t\ge0\) and every \(x\ge0\),
\begin{equation}\label{eq:generalized_lorden_tail}
\mathbb P\{W_t>x\}
\le
L_G\bigl(1-\Phi(x)\bigr)
+
\frac{1}{m_G}
\int_x^\infty \bigl(1-\Phi(u)\bigr)\,du.
\end{equation}
Consequently, for every \(0<\ell\le k-1\),
\begin{equation}\label{eq:generalized_lorden_moments}
\sup_{t\ge0}\mathbb E W_t^\ell
\le
L_G\,\mathbb E\zeta^\ell
+
\frac{\mathbb E\zeta^{\ell+1}}{(\ell+1)\,m_G}.
\end{equation}
In particular,
\begin{equation}\label{eq:generalized_lorden_firstmoment}
\sup_{t\ge0}\mathbb E W_t
\le
L_G\,\mathbb E\zeta
+
\frac{\mathbb E\zeta^{2}}{2\,m_G}.
\end{equation}
\end{theorem}

\begin{proof}
Fix \(t\ge0\) and \(x\ge0\). For \(n\ge0\), let
\[
\mathcal F_n:=\sigma(\zeta_1,\ldots,\zeta_n,\theta_1,\ldots,\theta_n).
\]
Then \(S_n\) is \(\mathcal F_n\)-measurable.

Since
\[
W_t=S_{N_t+1}-t,
\]
one has the disjoint decomposition
\[
\{W_t>x\}
=
\bigsqcup_{n=0}^\infty
\{S_n\le t,\ \xi_{n+1}>t+x-S_n\}.
\]
Indeed, on \(\{S_n\le t,\ \xi_{n+1}>t+x-S_n\}\), since \(x\ge0\),
\[
S_{n+1}=S_n+\xi_{n+1}>t+x\ge t,
\]
so \(S_n\le t<S_{n+1}\), hence \(N_t=n\).

Using Tonelli's theorem,
\begin{align*}
\mathbb P\{W_t>x\}
&=
\sum_{n=0}^\infty
\mathbb E\!\left[
\mathbf 1_{\{S_n\le t\}}
\mathbf 1_{\{\xi_{n+1}>t+x-S_n\}}
\right] \\
&=
\sum_{n=0}^\infty
\mathbb E\!\left[
\mathbf 1_{\{S_n\le t\}}
\mathbb P(\xi_{n+1}>t+x-S_n\mid\mathcal F_n)
\right].
\end{align*}
By Assumption~\ref{ass:grp:min},
\[
\xi_{n+1}\le \zeta_{n+1}
\qquad \text{almost surely},
\]
so
\[
\{\xi_{n+1}>y\}\subseteq \{\zeta_{n+1}>y\},
\qquad y\ge0.
\]
Since \(\zeta_{n+1}\) is independent of \(\mathcal F_n\),
\begin{align*}
\mathbb P\{W_t>x\}
&\le
\sum_{n=0}^\infty
\mathbb E\!\left[
\mathbf 1_{\{S_n\le t\}}
\mathbb P(\zeta_{n+1}>t+x-S_n)
\right] \\
&=
\sum_{n=0}^\infty
\mathbb E\!\left[
\mathbf 1_{\{S_n\le t\}}
(1-\Phi(t+x-S_n))
\right] \\
&=
\int_{[0,t]} (1-\Phi(t+x-s))\,U_\xi(ds).
\end{align*}
By Condition~\eqref{cond:RD},
\[
\mathbb P\{W_t>x\}
\le
\int_{[0,t]} (1-\Phi(t+x-s))\,U_G(ds).
\]

Define
\[
g_x(u):=1-\Phi(x+u),
\qquad u\ge0.
\]
Since \(\Phi\) is a distribution function and \(\mathbb E\zeta<\infty\), the function \(g_x\) is nonnegative, nonincreasing, right-continuous, and integrable. Moreover,
\[
g_x(0)=1-\Phi(x),
\qquad
\int_0^\infty g_x(u)\,du
=
\int_x^\infty (1-\Phi(v))\,dv.
\]
Applying Lemma~\ref{lem:renewal_measure_convolution} to \(g_x\), we obtain
\[
\mathbb P\{W_t>x\}
\le
L_G(1-\Phi(x))
+
\frac{1}{m_G}\int_x^\infty (1-\Phi(v))\,dv,
\]
which proves \eqref{eq:generalized_lorden_tail}.

For \(0<\ell\le k-1\), the tail-integral identity and \eqref{eq:generalized_lorden_tail} yield
\begin{align*}
\mathbb E W_t^\ell
&=
\ell\int_0^\infty x^{\ell-1}\mathbb P\{W_t>x\}\,dx \\
&\le
L_G\,
\ell\int_0^\infty x^{\ell-1}(1-\Phi(x))\,dx
+
\frac{\ell}{m_G}
\int_0^\infty x^{\ell-1}\int_x^\infty (1-\Phi(u))\,du\,dx.
\end{align*}
By Tonelli's theorem,
\[
\frac{\ell}{m_G}
\int_0^\infty x^{\ell-1}\int_x^\infty (1-\Phi(u))\,du\,dx
=
\frac{1}{m_G}
\int_0^\infty u^\ell(1-\Phi(u))\,du.
\]
Hence
\[
\mathbb E W_t^\ell
\le
L_G\,\mathbb E\zeta^\ell
+
\frac{\mathbb E\zeta^{\ell+1}}{(\ell+1)\,m_G},
\]
which proves \eqref{eq:generalized_lorden_moments}. The case \(\ell=1\) gives \eqref{eq:generalized_lorden_firstmoment}.
\end{proof}

\begin{remark}[Scope of the Lorden-type bound]
Theorem~\ref{thm:generalized_lorden} concerns the forward recurrence time \(W_t\), not the age process \(B_t\). Even in the classical setting, the coincidence of the limiting laws of \(B_t\) and \(W_t\) does not imply a uniform finite-time Lorden-type bound for \(B_t\).
\end{remark}

\begin{remark}[What must be checked in applications]
Condition~\eqref{cond:RD} is the additional ingredient that turns one-dimensional comparison of inter-renewal laws into a uniform upper bound for \(W_t\). It is not a consequence of the marginal inequalities
\[
G(x)\ge F_i(x)\ge \Phi(x),
\qquad x\ge0,
\]
and must therefore be verified separately, or replaced by another structural condition yielding the same conclusion.
\end{remark}

\subsection{Consequences for coupling and convergence rates}

Set
\[
C_{\mathrm{GL}}
:=
L_G\,\mathbb E\zeta
+
\frac{\mathbb E\zeta^{2}}{2\,m_G},
\qquad
m_G=\mathbb E\eta,
\qquad
L_G=\frac{\mathbb E\eta^2}{(\mathbb E\eta)^2}.
\]
Then \eqref{eq:generalized_lorden_firstmoment} gives
\[
\sup_{t\ge0}\mathbb E W_t\le C_{\mathrm{GL}}.
\]
\rev{Thus Theorem~\ref{thm:generalized_lorden} supplies the first-moment input used in the coupling constructions of Sections~\ref{sec:stateOftheArt}--\ref{section:classicalresults}. The next corollary records the corresponding total-variation estimate only conditionally, once the remaining coupling ingredients are available.}

\begin{corollary}[Convergence-rate bound from the generalised Lorden estimate]
\label{cor:generalized_lorden_coupling}
Let \(X=(X_t)_{t\ge0}\) be a process on a measurable state space \(E\). Fix \(x\in E\), and suppose that one can construct on a common filtered probability space
\[
(\Omega,\mathcal F,(\mathcal F_t)_{t\ge0},\mathbb P)
\]
\begin{itemize}
    \item a version \(X^x=(X_t^x)_{t\ge0}\) of the process started from \(x\),
    \item a stationary version \(\widetilde X=(\widetilde X_t)_{t\ge0}\) of the same process, whose one-time marginal law is \(\pi\),
    \item an increasing sequence of almost surely finite stopping times
    \[
    0=\Gamma_0<\Gamma_1<\Gamma_2<\cdots,
    \]
    and events \(A_n\in\mathcal F_{\Gamma_n}\), \(n\ge1\),
\end{itemize}
such that:
\begin{enumerate}
    \item
    \[
    \mathbb E\Gamma_1\le C_{\mathrm{GL}};
    \]
    \item for every \(n\ge1\),
    \[
    \mathbb E[\Gamma_{n+1}-\Gamma_n\mid\mathcal F_{\Gamma_n}]\le M_x
    \qquad \text{a.s.},
    \]
    for some \(M_x<\infty\);
    \item for every \(n\ge1\),
    \[
    \mathbb P(A_n\mid\mathcal F_{\Gamma_{n-1}})
    \ge q_x
    \qquad \text{a.s. on } \bigcap_{k=1}^{n-1}A_k^c,
    \]
    for some \(q_x\in(0,1]\);
    \item on \(A_n\),
    \[
    X_t^x=\widetilde X_t
    \qquad \text{for all } t\ge \Gamma_n
    \quad \text{a.s.}
    \]
\end{enumerate}
Then, with
\[
N:=\inf\{n\ge1:\ A_n \text{ occurs}\},
\qquad
\tau_x:=\Gamma_N,
\]
one has
\[
\mathbb E\tau_x
\le
C_{\mathrm{GL}}+M_x(q_x^{-1}-1),
\]
and consequently
\[
\|\Law(X_t^x)-\pi\|_{\mathrm{TV}}
\le
\frac{C_{\mathrm{GL}}+M_x(q_x^{-1}-1)}{t},
\qquad t>0.
\]
\end{corollary}

\begin{proof}
The argument is the same as in Proposition~\ref{prop:trial_coupling_classical}, with \((T_n,a_b,m,q_b)\) replaced by \((\Gamma_n,C_{\mathrm{GL}},M_x,q_x)\). This gives
\[
\mathbb P(N>n)\le (1-q_x)^n,
\qquad n\ge0,
\]
and hence
\[
\mathbb E\tau_x
=
\mathbb E\Gamma_N
\le
C_{\mathrm{GL}}+M_x(q_x^{-1}-1).
\]
The final total-variation estimate follows from Proposition~\ref{prop:coupling_stationary} with \(\varphi(u)=u\).
\end{proof}

\begin{remark}[How Theorem~\ref{thm:generalized_lorden} enters Corollary~\ref{cor:generalized_lorden_coupling}]
The condition
\[
\mathbb E\Gamma_1\le C_{\mathrm{GL}}
\]
is the point at which the Lorden-type estimate enters the coupling construction. In concrete models, \(\Gamma_1\) is typically chosen as an initial trial delay controlled by the forward recurrence time, and the above inequality is then verified either by Theorem~\ref{thm:generalized_lorden} or by a model-specific comparison argument.
\end{remark}

\begin{remark}[Transfer through the age-process kernel representation]
\label{rem:generalized_lorden_transfer}
If the coupling is first constructed for the associated age process and the one-time law of \(X_t\) admits a time-independent kernel representation through that age process, then Lemma~\ref{lem:sxodim} transfers the total-variation bound directly to \(X_t\). In particular, if
\[
\mathcal L(X_t)(A)=\int_{[0,\infty)}K(s,A)\,\mathcal L(B_t^b)(ds),
\qquad A\in\mathcal B(E),
\]
for the same kernel \(K\) as in Condition~\eqref{eq:ConditionK}, and if
\[
\pi_X(A):=\int_{[0,\infty)}K(s,A)\,\pi_B(ds),
\qquad A\in\mathcal B(E),
\]
has been verified to be the one-time marginal law of a stationary version of \(X\), then
\[
\|\mathcal L(X_t)-\pi_X\|_{\mathrm{TV}}
\le
\|\mathcal L(B_t^b)-\pi_B\|_{\mathrm{TV}},
\]
so the convergence-rate estimate for \(X_t\) follows from the corresponding estimate for the age process.
\end{remark}

\section{Illustrative examples and comparison with classical bounds}
\label{sec:application}

\rev{This section records four benchmark calculations: exponential, mixed, Markov-modulated, and Pareto. In the i.i.d.\ cases, Condition~\eqref{cond:RD} holds automatically when the comparison variables are chosen to have the same law as the inter-renewal times. In the Markov-modulated benchmark, only the comparison assumptions and the explicit constant are verified; the final convergence consequence remains conditional on \eqref{cond:RD}. The comparative tables are deferred to Appendix~\ref{app:benchmark_tables}.}

\subsection{Benchmark light-tailed model: exponential inter-renewal times}
\label{subsec:benchmark_exponential}

\rev{We begin with the exponential benchmark, where both the forward recurrence time and the total-variation distance to stationarity are explicit. This calibrates the universal constant against an exactly solvable model.}

Assume throughout this subsection that
\[
\xi_1,\xi_2,\ldots \stackrel{\mathrm{i.i.d.}}{\sim} \mathrm{Exp}(\lambda),
\qquad \lambda>0,
\]
and let \(N=(N_t)_{t\ge0}\) be the associated renewal process. Since the exponential law is memoryless, \(N\) is the standard Poisson process of rate \(\lambda\).

\begin{proposition}[Exact forward recurrence and exact convergence rate in the exponential case]
\label{prop:benchmark_exponential_exact}
Let \(W_t\) be the forward recurrence time and let \(B_t^0\) be the age process started from the renewal state \(0\). Let \(\pi_B\) denote the stationary law of the age process. Then the following hold.
\begin{enumerate}
    \item For every \(t\ge0\),
    \[
    W_t \overset{D}{=} \mathrm{Exp}(\lambda),
    \qquad
    \mathbb E W_t=\frac1\lambda.
    \]

    \item The stationary age law is
    \[
    \pi_B(dy)=\lambda e^{-\lambda y}\,dy,
    \qquad y\ge0.
    \]

    \item For every \(t>0\),
    \[
    \Law(B_t^0)(dy)
    =
    \lambda e^{-\lambda y}\mathbf 1_{[0,t)}(y)\,dy
    +
    e^{-\lambda t}\delta_t(dy).
    \]

    \item Consequently,
    \[
    \|\Law(B_t^0)-\pi_B\|_{\mathrm{TV}}
    =
    e^{-\lambda t},
    \qquad t\ge0.
    \]
\end{enumerate}
\end{proposition}

\begin{proof}
For every \(x\ge0\),
\[
\mathbb P\{W_t>x\}
=
\mathbb P\{N_{t+x}-N_t=0\}
=
e^{-\lambda x},
\]
because the Poisson process has stationary independent increments. This proves (1).

The stationary age law of a Poisson renewal process is exponential with parameter \(\lambda\), and also follows from the equilibrium formula of Section~\ref{subsec:stateOftheArtLordenClassic}. This gives (2).

To prove (3), fix \(t>0\). For \(0\le y<t\),
\[
\mathbb P\{B_t^0>y\}
=
\mathbb P\{N_t-N_{t-y}=0\}
=
e^{-\lambda y},
\]
again by stationary independent increments. Thus \(B_t^0\) has density
\[
\lambda e^{-\lambda y},
\qquad 0\le y<t.
\]
The remaining probability mass is
\[
\mathbb P\{B_t^0=t\}=\mathbb P\{N_t=0\}=e^{-\lambda t},
\]
and
\[
\int_0^t \lambda e^{-\lambda y}\,dy + e^{-\lambda t}
=
(1-e^{-\lambda t})+e^{-\lambda t}
=
1,
\]
so the displayed measure is indeed a probability law.

For (4), the laws \(\Law(B_t^0)\) and \(\pi_B\) have the same density on \([0,t)\), while \(\Law(B_t^0)\) has an atom of mass \(e^{-\lambda t}\) at \(t\), and \(\pi_B\) assigns mass
\[
\pi_B((t,\infty))=e^{-\lambda t}
\]
to the tail \((t,\infty)\). Hence
\[
\|\Law(B_t^0)-\pi_B\|_{\mathrm{TV}}=e^{-\lambda t}.
\]
\end{proof}

\begin{proposition}[A sharp coalescent coupling in the exponential benchmark]
\label{prop:benchmark_exponential_coupling}
In the setting of Proposition~\ref{prop:benchmark_exponential_exact}, there exists a coupling of \(B^0=(B_t^0)_{t\ge0}\) with a stationary version \(\widetilde B=(\widetilde B_t)_{t\ge0}\) such that the coalescence time \(\tau\) satisfies
\[
\tau \overset{D}{=} \mathrm{Exp}(\lambda).
\]
Consequently,
\[
\|\Law(B_t^0)-\pi_B\|_{\mathrm{TV}}
\le
\mathbb P\{\tau>t\}
=
e^{-\lambda t},
\qquad t\ge0,
\]
and this bound is exact by Proposition~\ref{prop:benchmark_exponential_exact}.
\end{proposition}

\begin{proof}
Let \(A\) and \(E\) be independent random variables such that
\[
A\sim \mathrm{Exp}(\lambda),
\qquad
E\sim \mathrm{Exp}(\lambda).
\]
Let \(\Gamma=(\Gamma_n)_{n\ge1}\) be an independent sequence of arrival times of a Poisson process of rate \(\lambda\) on \((0,\infty)\), with the convention \(\Gamma_0:=0\).

Construct \(B^0\) by placing its renewal epochs at
\[
E+\Gamma_n,\qquad n\ge0.
\]
Construct \(\widetilde B\) by placing the last renewal before \(0\) at time \(-A\), the first renewal after \(0\) at time \(E\), and all subsequent renewal epochs at
\[
E+\Gamma_n,\qquad n\ge1.
\]
Then \(\widetilde B_0=A\sim\pi_B\), and the two processes use the same renewal epochs from time \(E\) onward. Hence
\[
B_t^0=\widetilde B_t
\qquad \text{for all } t\ge E,
\]
so one may take \(\tau:=E\sim \mathrm{Exp}(\lambda)\).
\end{proof}

\begin{remark}[Calibration of the universal Lorden input]
\label{rem:benchmark_exponential_calibration}
In this benchmark,
\[
\sup_{t\ge0}\mathbb E W_t=\frac1\lambda,
\]
whereas the universal constant of Section~\ref{sec:mainresultLordenGener} becomes
\[
C_{\mathrm{GL}}^{\mathrm{exp}}
=
\frac{3}{\lambda}.
\]
Thus the generalised Lorden input has the correct scale \(1/\lambda\), but is conservative by the factor \(3\).
\end{remark}

\begin{remark}[Generic versus refined use of the coupling method]
\label{rem:benchmark_exponential_generic_vs_refined}
If one uses only the first moment of the coalescence time \(\tau\sim\mathrm{Exp}(\lambda)\), then Proposition~\ref{prop:coupling_stationary} gives
\[
\|\Law(B_t^0)-\pi_B\|_{\mathrm{TV}}
\le
\frac{1}{\lambda t},
\qquad t>0.
\]
Using instead the exact tail of \(\tau\), the explicit coupling of Proposition~\ref{prop:benchmark_exponential_coupling} recovers the exact rate
\[
\|\Law(B_t^0)-\pi_B\|_{\mathrm{TV}}
=
e^{-\lambda t}.
\]
This benchmark therefore separates the universal first-moment input from sharper model-specific rates.
\end{remark}

\subsection{Benchmark mixed law: an atom plus an exponential tail}
\label{subsec:benchmark_mixed}

\rev{We next consider a genuinely mixed inter-renewal law. The purpose is to show that the measure-valued formulation of Section~\ref{sec:mainresultLordenGener} remains explicit in the presence of an atom.}

Assume throughout this subsection that the inter-renewal times are i.i.d.\ with distribution function
\begin{equation}\label{eq:mixed_benchmark_distribution}
F(x)=
\begin{cases}
0, & x<d,\\[1mm]
p, & x=d,\\[1mm]
1-(1-p)e^{-\lambda(x-d)}, & x>d,
\end{cases}
\qquad
d>0,\quad p\in(0,1),\quad \lambda>0.
\end{equation}
Equivalently,
\[
\mathbb P\{\xi=d\}=p,
\qquad
\mathbb P\{\xi>d+x\}=(1-p)e^{-\lambda x},
\qquad x\ge0.
\]

\begin{proposition}[Mixed-law benchmark: explicit generalised intensity and Lorden input]
\label{prop:benchmark_mixed}
Let \(\xi\) have distribution \eqref{eq:mixed_benchmark_distribution}. Then the following hold.
\begin{enumerate}
    \item The generalised intensity measure of \(\xi\) is
    \[
    \lambda_\xi(ds)=\alpha\,\delta_d(ds)+\lambda\,\mathbf 1_{(d,\infty)}(s)\,ds,
    \qquad
    \alpha:=-\log(1-p).
    \]

    \item The first two moments are
    \[
    \mathbb E\xi
    =
    d+\frac{1-p}{\lambda},
    \]
    and
    \[
    \mathbb E\xi^2
    =
    d^2+\frac{2(1-p)d}{\lambda}+\frac{2(1-p)}{\lambda^2}.
    \]

    \item The equilibrium mean forward recurrence time is
    \[
    \mathbb E W_\infty
    =
    \frac{\mathbb E\xi^2}{2\,\mathbb E\xi}
    =
    \frac{1}{2}\,
    \frac{d^2+\frac{2(1-p)d}{\lambda}+\frac{2(1-p)}{\lambda^2}}
         {d+\frac{1-p}{\lambda}}.
    \]

    \item If the comparison variables are chosen so that
    \[
    \eta\overset{D}{=}\xi\overset{D}{=}\zeta,
    \]
    then
    \[
    m_G=\mathbb E\eta=\mathbb E\xi,
    \qquad
    L_G=\frac{\mathbb E\eta^2}{(\mathbb E\eta)^2}
       =\frac{\mathbb E\xi^2}{(\mathbb E\xi)^2},
    \]
    and therefore
    \[
    C_{\mathrm{GL}}^{\mathrm{mix}}
    =
    \frac{3}{2}\frac{\mathbb E\xi^2}{\mathbb E\xi}.
    \]
    Equivalently,
    \[
    C_{\mathrm{GL}}^{\mathrm{mix}}
    =
    \frac{3}{2}\,
    \frac{d^2+\frac{2(1-p)d}{\lambda}+\frac{2(1-p)}{\lambda^2}}
         {d+\frac{1-p}{\lambda}}.
    \]
\end{enumerate}
\end{proposition}

\begin{proof}
For \(x<d\), one has \(\overline F(x)=1\). At the atom \(d\),
\[
\overline F(d-)=1,
\qquad
\overline F(d)=1-p,
\]
so the atomic contribution to the generalised intensity measure is
\[
-\log\frac{\overline F(d)}{\overline F(d-)}
=
-\log(1-p)
=
\alpha.
\]
For \(x>d\),
\[
\overline F(x)=(1-p)e^{-\lambda(x-d)},
\]
so the absolutely continuous hazard on \((d,\infty)\) is \(\lambda\). This proves (1).

For the first moment,
\[
\mathbb E\xi
=
p\,d+(1-p)\,\mathbb E(d+E)
=
d+\frac{1-p}{\lambda},
\]
where \(E\sim\mathrm{Exp}(\lambda)\). Likewise,
\[
\mathbb E\xi^2
=
p\,d^2+(1-p)\,\mathbb E(d+E)^2
=
d^2+\frac{2(1-p)d}{\lambda}+\frac{2(1-p)}{\lambda^2},
\]
which proves (2).

Since the benchmark is i.i.d.,
\[
\mathbb E W_\infty=\frac{\mathbb E\xi^2}{2\,\mathbb E\xi},
\]
and substitution of the moments gives (3).

For (4), if \(\eta\), \(\xi\), and \(\zeta\) have the same law, then
\[
m_G=\mathbb E\xi,
\qquad
L_G=\frac{\mathbb E\xi^2}{(\mathbb E\xi)^2},
\]
and therefore
\[
C_{\mathrm{GL}}^{\mathrm{mix}}
=
L_G\,\mathbb E\xi+\frac{\mathbb E\xi^2}{2\,\mathbb E\xi}
=
\frac{3}{2}\frac{\mathbb E\xi^2}{\mathbb E\xi}.
\]
Substituting the explicit moments completes the proof.
\end{proof}

\begin{remark}[Calibration against the classical i.i.d.\ Lorden scale]
\label{rem:benchmark_mixed_calibration}
For the mixed law \eqref{eq:mixed_benchmark_distribution},
\[
\frac{\mathbb E\xi^2}{\mathbb E\xi}
\]
is the classical i.i.d.\ Lorden scale, while the present construction gives
\[
C_{\mathrm{GL}}^{\mathrm{mix}}
=
\frac{3}{2}\frac{\mathbb E\xi^2}{\mathbb E\xi}.
\]
Thus the correct renewal-theoretic scale is preserved, up to the same universal factor \(3/2\) as in the exponential benchmark.
\end{remark}

\begin{remark}[Why the mixed-law benchmark matters]
\label{rem:benchmark_mixed_why}
The law \eqref{eq:mixed_benchmark_distribution} lies outside the purely absolutely continuous framework of Section~\ref{section:classicalresults}. Nevertheless, the generalised intensity measure and the resulting Lorden-type input remain explicit.
\end{remark}

\begin{remark}[Consequence for convergence-rate bounds]
\label{rem:benchmark_mixed_rate}
Once a coupling construction is available and the trial parameters \(m_x\) and \(q_x\) are identified, Corollary~\ref{cor:generalized_lorden_coupling} yields
\[
\|\Law(X_t^x)-\pi\|_{\mathrm{TV}}
\le
\frac{C_{\mathrm{GL}}^{\mathrm{mix}}+m_x(q_x^{-1}-1)}{t},
\qquad t>0.
\]
In this i.i.d.\ benchmark, if \(\eta\overset{D}{=}\xi\overset{D}{=}\zeta\), then \(U_\xi=U_G\), so Condition~\eqref{cond:RD} holds automatically.
\end{remark}

\subsection{Dependent light-tailed benchmark for the comparison framework}
\label{subsec:benchmark_mme}

\rev{We now turn to a dependent light-tailed benchmark. The point is to verify the comparison assumptions and the explicit constant $C_{\mathrm{GL}}^{\mathrm{MME}}$ in closed form under dependence and heterogeneity.}

Let \(J=(J_i)_{i\ge1}\) be a time-homogeneous Markov chain on the finite state space
\[
\mathcal S=\{1,\dots,m\},
\]
with arbitrary initial distribution and transition matrix \(P\). Let
\[
0<\underline\lambda\le \lambda_j\le \overline\lambda<\infty,
\qquad j\in\mathcal S,
\]
and fix
\[
0<\phi<\underline\lambda.
\]

Let \(\zeta_1,\zeta_2,\ldots\) be i.i.d.\ random variables with
\[
\zeta_i\sim \mathrm{Exp}(\phi),
\qquad i\ge1,
\]
independent of \(J\). Conditionally on \(J\), let \(\theta_1,\theta_2,\ldots\) be independent random variables such that
\[
\theta_i\mid \{J_i=j\}\sim \mathrm{Exp}(\lambda_j-\phi),
\qquad j\in\mathcal S,
\]
and assume that the family \(\{\theta_i\}_{i\ge1}\) is independent of the sequence \(\{\zeta_i\}_{i\ge1}\). Define
\[
\xi_i:=\min\{\zeta_i,\theta_i\},
\qquad i\ge1.
\]

Conditionally on \(\{J_i=j\}\), the variable \(\xi_i\) is exponential with parameter \(\lambda_j\). Thus the sequence \((\xi_i)\) is, in general, dependent and non-identically distributed through the modulating chain \(J\).

\begin{proposition}[Verification of the Section~\ref{sec:mainresultLordenGener} assumptions]
\label{prop:benchmark_mme_assumptions}
In the Markov-modulated exponential benchmark described above, Assumptions~\ref{ass:grp:min}--\ref{ass:grp:moment} are satisfied with
\[
\varphi(ds)=\phi\,ds,
\qquad
Q(ds)=\overline\lambda\,ds.
\]
Moreover,
\[
\Phi(x)=1-e^{-\phi x},
\qquad
G(x)=1-e^{-\overline\lambda x},
\qquad x\ge0.
\]
\end{proposition}

\begin{proof}
Assumption~\ref{ass:grp:min} holds by construction.

For each \(i\ge1\), let
\[
p_i(j):=\mathbb P\{J_i=j\},
\qquad j\in\mathcal S.
\]
Then
\[
\mathbb P\{\theta_i>x\}
=
\sum_{j=1}^m p_i(j)e^{-(\lambda_j-\phi)x},
\qquad x\ge0.
\]
Hence \(\theta_i\) is absolutely continuous on \((0,\infty)\), with hazard rate
\[
h_i(x)
=
\frac{\sum_{j=1}^m p_i(j)(\lambda_j-\phi)e^{-(\lambda_j-\phi)x}}
     {\sum_{j=1}^m p_i(j)e^{-(\lambda_j-\phi)x}},
\qquad x\ge0.
\]
This is a convex combination of the values \(\lambda_j-\phi\), so
\[
0\le h_i(x)\le \overline\lambda-\phi,
\qquad x\ge0.
\]
Thus
\[
\mu_i(ds)=h_i(s)\,ds,
\qquad
\varphi+\mu_i\le \overline\lambda\,ds=Q(ds),
\]
which is Assumption~\ref{ass:grp:upper}.

Assumption~\ref{ass:grp:proper} holds since
\[
\varphi((0,\infty))=\int_0^\infty \phi\,ds=\infty.
\]
Assumption~\ref{ass:grp:positivity} holds with \(T=0\), because
\[
\varphi^{\mathrm{ac}}(s)\equiv\phi>0
\qquad \text{for all } s\ge0.
\]
Finally,
\[
\int_0^\infty x^{k-1}\exp\!\bigl(-\varphi((0,x])\bigr)\,dx
=
\int_0^\infty x^{k-1}e^{-\phi x}\,dx
<\infty,
\]
for every \(k\ge2\), so Assumption~\ref{ass:grp:moment} also holds. The formulas for \(\Phi\) and \(G\) follow from Proposition~\ref{prop:reconstruction_generalised_intensity}.
\end{proof}

\begin{corollary}[Explicit generalised Lorden constant in the Markov-modulated benchmark]
\label{cor:benchmark_mme_constant}
In the setting of Proposition~\ref{prop:benchmark_mme_assumptions},
\[
m_G=\mathbb E\eta=\frac{1}{\overline\lambda},
\qquad
L_G=\frac{\mathbb E\eta^2}{(\mathbb E\eta)^2}=2,
\]
and therefore
\[
C_{\mathrm{GL}}^{\mathrm{MME}}
=
L_G\,\mathbb E\zeta
+
\frac{\mathbb E\zeta^2}{2\,m_G}
=
\frac{2}{\phi}+\frac{\overline\lambda}{\phi^2}.
\]
\end{corollary}

\begin{proof}
Since \(\eta\sim \mathrm{Exp}(\overline\lambda)\),
\[
\mathbb E\eta=\frac{1}{\overline\lambda},
\qquad
\mathbb E\eta^2=\frac{2}{\overline\lambda^2},
\]
so
\[
L_G=2.
\]
Likewise, since \(\zeta\sim \mathrm{Exp}(\phi)\),
\[
\mathbb E\zeta=\frac{1}{\phi},
\qquad
\mathbb E\zeta^2=\frac{2}{\phi^2}.
\]
Substituting into the definition of \(C_{\mathrm{GL}}\) yields the claim.
\end{proof}

\begin{proposition}[Renewal-function domination in the Markov-modulated benchmark]
\label{prop:benchmark_mme_renewal_function}
In the setting of Subsection~\ref{subsec:benchmark_mme}, one can realise the benchmark on an enlarged probability space in such a way that
\[
\eta_i\le \xi_i
\qquad \text{almost surely for every } i\ge1,
\]
where the variables \(\eta_i\) are i.i.d.\ with law \(G=\mathrm{Exp}(\overline\lambda)\). Consequently,
\[
S_n^\eta\le S_n^\xi
\qquad \text{almost surely for every } n\ge1,
\]
and therefore
\[
U_\xi([0,t])\le U_G([0,t]),
\qquad t\ge0.
\]
\end{proposition}

\begin{proof}
Let \(E_1,E_2,\ldots\stackrel{\mathrm{i.i.d.}}{\sim}\mathrm{Exp}(1)\) be independent of \(J\), and define
\[
\eta_i:=\frac{E_i}{\overline\lambda},
\qquad
\widehat\xi_i:=\frac{E_i}{\lambda_{J_i}},
\qquad i\ge1.
\]
Then \(\eta_i\sim \mathrm{Exp}(\overline\lambda)\), and conditionally on \(\{J_i=j\}\),
\[
\widehat\xi_i\mid \{J_i=j\}\sim \mathrm{Exp}(\lambda_j).
\]
Since \(\lambda_{J_i}\le \overline\lambda\),
\[
\eta_i=\frac{E_i}{\overline\lambda}\le \frac{E_i}{\lambda_{J_i}}=\widehat\xi_i
\qquad \text{almost surely}.
\]
Hence
\[
S_n^\eta:=\sum_{k=1}^n \eta_k
\le
\sum_{k=1}^n \widehat\xi_k=:\widehat S_n^\xi
\qquad \text{almost surely}.
\]
Therefore
\[
\{\widehat S_n^\xi\le t\}\subseteq \{S_n^\eta\le t\},
\qquad t\ge0,
\]
and so
\[
\mathbb P\{S_n^\xi\le t\}
=
\mathbb P\{\widehat S_n^\xi\le t\}
\le
\mathbb P\{S_n^\eta\le t\}.
\]
Summing over \(n\ge0\) yields
\[
U_\xi([0,t])\le U_G([0,t]),
\qquad t\ge0.
\]
\end{proof}

\begin{remark}[What this does and does not prove]
The previous proposition yields only cumulative renewal-function domination,
\[
U_\xi([0,t])\le U_G([0,t]),
\qquad t\ge0.
\]
This is weaker than the setwise domination condition
\[
U_\xi(A)\le U_G(A),
\qquad A\in\mathcal B([0,\infty)),
\]
required in Condition~\eqref{cond:RD}. Accordingly, the explicit constant \(C_{\mathrm{GL}}^{\mathrm{MME}}\) is fully verified in this benchmark, whereas the final convergence-rate theorem remains conditional on Condition~\eqref{cond:RD}.
\end{remark}

\begin{remark}[Comparison with sharper model-specific methods]
\label{rem:benchmark_mme_comparison}
\rev{In the finite-state Markov-modulated setting, matrix-analytic or spectral methods may yield sharper rates, but they depend on the full transition structure of the modulating chain. Here the benchmark is used only to display the comparison constant in closed form.}
\end{remark}

\begin{remark}[Correct order and degeneration to the classical benchmark]
\label{rem:benchmark_mme_scale}
The constant \(C_{\mathrm{GL}}^{\mathrm{MME}}\) has the natural light-tailed scale \(O(\phi^{-1})\) when \(\overline\lambda\) is comparable to \(\phi\). If \(\lambda_j\equiv\lambda\) for all \(j\), then choosing \(\phi\uparrow\lambda\) gives
\[
C_{\mathrm{GL}}^{\mathrm{MME}}\to \frac{3}{\lambda},
\]
which coincides with the universal constant in the exponential benchmark.
\end{remark}

\begin{remark}[Why this benchmark matters]
\label{rem:benchmark_mme_why}
\rev{This benchmark shows that dependence does not prevent closed-form verification of the comparison assumptions and of \(C_{\mathrm{GL}}^{\mathrm{MME}}\). The full convergence-rate conclusion remains conditional on Condition~\eqref{cond:RD}.}
\end{remark}

\subsection{Benchmark heavy-tailed model: Pareto inter-renewal times}
\label{subsec:benchmark_pareto}

\rev{We finally consider a heavy-tailed i.i.d.\ benchmark. The focus is the Lorden input and the associated moment threshold, rather than an exact rate.}

Assume throughout this subsection that the inter-renewal times are i.i.d.\ with distribution function
\begin{equation}\label{eq:pareto_benchmark_distribution}
F(x)=1-\left(1+\frac{x}{b}\right)^{-\alpha},
\qquad x\ge0,
\qquad b>0,\quad \alpha>2.
\end{equation}
Equivalently,
\[
\mathbb P\{\xi>x\}
=
\left(1+\frac{x}{b}\right)^{-\alpha},
\qquad x\ge0.
\]

\begin{proposition}[Pareto benchmark: explicit generalised intensity and moment threshold]
\label{prop:benchmark_pareto}
Let \(\xi\) have distribution \eqref{eq:pareto_benchmark_distribution}. Then the following hold.
\begin{enumerate}
    \item The generalised intensity measure of \(\xi\) is purely absolutely continuous and is given by
    \[
    \lambda_\xi(ds)=\frac{\alpha}{b+s}\,ds.
    \]

    \item For every \(r\in(0,\alpha)\),
    \[
    \mathbb E\xi^r
    =
    b^r\frac{\Gamma(r+1)\Gamma(\alpha-r)}{\Gamma(\alpha)}.
    \]
    In particular,
    \[
    \mathbb E\xi=\frac{b}{\alpha-1},
    \qquad
    \mathbb E\xi^2=\frac{2b^2}{(\alpha-1)(\alpha-2)}.
    \]

    \item For every integer \(k\) such that \(2\le k<\alpha\), Assumption~\ref{ass:grp:moment} is satisfied with
    \[
    \varphi(ds)=\frac{\alpha}{b+s}\,ds.
    \]

    \item If the comparison variables are chosen so that
    \[
    \eta\overset{D}{=}\xi\overset{D}{=}\zeta,
    \]
    then
    \[
    m_G=\mathbb E\eta=\frac{b}{\alpha-1},
    \qquad
    L_G=\frac{\mathbb E\eta^2}{(\mathbb E\eta)^2}
       =\frac{2(\alpha-1)}{\alpha-2},
    \]
    and hence
    \[
    C_{\mathrm{GL}}^{\mathrm{Par}}
    =
    \frac{3b}{\alpha-2}.
    \]
\end{enumerate}
\end{proposition}

\begin{proof}
Since
\[
\mathbb P\{\xi>x\}
=
\left(1+\frac{x}{b}\right)^{-\alpha}
=
\exp\!\left(-\int_0^x \frac{\alpha}{b+s}\,ds\right),
\qquad x\ge0,
\]
Proposition~\ref{prop:reconstruction_generalised_intensity} gives
\[
\lambda_\xi(ds)=\frac{\alpha}{b+s}\,ds,
\]
which proves (1).

For \(r\in(0,\alpha)\), the density of \(\xi\) is
\[
f(x)=\frac{\alpha}{b}\left(1+\frac{x}{b}\right)^{-\alpha-1},
\qquad x\ge0.
\]
Therefore,
\[
\mathbb E\xi^r
=
\frac{\alpha}{b}\int_0^\infty x^r\left(1+\frac{x}{b}\right)^{-\alpha-1}\,dx.
\]
With the change of variables \(x=bu\),
\[
\mathbb E\xi^r
=
\alpha b^r\int_0^\infty u^r(1+u)^{-\alpha-1}\,du.
\]
By the standard Beta-integral formula,
\[
\int_0^\infty u^r(1+u)^{-\alpha-1}\,du
=
B(r+1,\alpha-r)
=
\frac{\Gamma(r+1)\Gamma(\alpha-r)}{\Gamma(\alpha+1)},
\]
and using \(\Gamma(\alpha+1)=\alpha\,\Gamma(\alpha)\), we obtain
\[
\mathbb E\xi^r
=
b^r\frac{\Gamma(r+1)\Gamma(\alpha-r)}{\Gamma(\alpha)}.
\]
The cases \(r=1\) and \(r=2\) give
\[
\mathbb E\xi=\frac{b}{\alpha-1},
\qquad
\mathbb E\xi^2=\frac{2b^2}{(\alpha-1)(\alpha-2)}.
\]

For (3), note that
\[
\exp\!\bigl(-\varphi((0,x])\bigr)
=
\left(1+\frac{x}{b}\right)^{-\alpha},
\]
so Assumption~\ref{ass:grp:moment} becomes
\[
\int_0^\infty x^{k-1}\left(1+\frac{x}{b}\right)^{-\alpha}\,dx<\infty.
\]
This holds exactly when \(k<\alpha\).

Finally, if \(\eta\), \(\xi\), and \(\zeta\) have the same law, then
\[
m_G=\mathbb E\xi=\frac{b}{\alpha-1},
\qquad
L_G=\frac{\mathbb E\xi^2}{(\mathbb E\xi)^2}
=
\frac{2(\alpha-1)}{\alpha-2}.
\]
Substituting into the definition of \(C_{\mathrm{GL}}\) gives
\[
C_{\mathrm{GL}}^{\mathrm{Par}}
=
L_G\,\mathbb E\xi+\frac{\mathbb E\xi^2}{2\,\mathbb E\xi}
=
\frac{3b}{\alpha-2}.
\]
\end{proof}

\begin{remark}[Correct threshold for the first Lorden input]
\label{rem:benchmark_pareto_threshold}
The benchmark constant
\[
C_{\mathrm{GL}}^{\mathrm{Par}}=\frac{3b}{\alpha-2}
\]
is finite if and only if \(\alpha>2\), which is exactly the natural second-moment threshold for the inter-renewal law.
\end{remark}

\begin{remark}[Higher moments and polynomial regimes]
\label{rem:benchmark_pareto_polynomial}
\rev{For every integer \(k\) with \(2\le k<\alpha\), Theorem~\ref{thm:generalized_lorden} yields uniform bounds on \(\mathbb E W_t^\ell\) for all \(0<\ell\le k-1\). If the remaining coupling terms admit moments of the same order, the resulting coupling bounds are polynomial.}
\end{remark}

\begin{remark}[Comparison with generator/Lyapunov approaches]
\label{rem:benchmark_pareto_lyapunov}
\rev{In the i.i.d.\ Pareto benchmark one may also pursue generator-based or Foster--Lyapunov arguments for the age process. Here the quantitative input is read off directly from the comparison intensity}
\[
\varphi(ds)=\frac{\alpha}{b+s}\,ds.
\]
\end{remark}

\begin{remark}[Use in convergence-rate estimates]
\label{rem:benchmark_pareto_rate}
Once a coupling construction is available and the trial parameters \(m_x\) and \(q_x\) are identified, Corollary~\ref{cor:generalized_lorden_coupling} yields
\[
\|\Law(X_t^x)-\pi\|_{\mathrm{TV}}
\le
\frac{C_{\mathrm{GL}}^{\mathrm{Par}}+m_x(q_x^{-1}-1)}{t},
\qquad t>0.
\]
If higher moments of the coupling time are available, Remark~\ref{rem:benchmark_pareto_polynomial} leads to higher-order polynomial rates.
\end{remark}

\subsection{Comparative summary and benchmark tables}
\label{subsec:benchmark_summary}

\rev{The benchmark comparison separates the explicit Lorden input from the subsequent coupling step. In the i.i.d.\ setting, Proposition~\ref{prop:benchmark_factor_three} shows that the universal constant is exactly \(3\,\mathbb E W_\infty\). For readability, the detailed comparative tables are collected in Appendix~\ref{app:benchmark_tables}.}

\begin{proposition}[Universal calibration in the i.i.d.\ benchmark]
\label{prop:benchmark_factor_three}
Suppose that the inter-renewal times are i.i.d.\ with common distribution \(F\), that \(\mathbb E\xi_1^2<\infty\), and that the comparison variables are chosen so that
\[
\eta\overset{D}{=}\xi_1\overset{D}{=}\zeta.
\]
Then
\[
C_{\mathrm{GL}}
=
\frac{3}{2}\frac{\mathbb E\xi_1^2}{\mathbb E\xi_1}
=
3\,\mathbb E W_\infty,
\]
where \(W_\infty\) denotes the equilibrium forward recurrence time of the classical renewal process.
\end{proposition}

\begin{proof}
In this case,
\[
m_G=\mathbb E\xi_1,
\qquad
L_G=\frac{\mathbb E\xi_1^2}{(\mathbb E\xi_1)^2},
\]
and therefore
\[
C_{\mathrm{GL}}
=
L_G\,\mathbb E\xi_1+\frac{\mathbb E\xi_1^2}{2\,\mathbb E\xi_1}
=
\frac{3}{2}\frac{\mathbb E\xi_1^2}{\mathbb E\xi_1}.
\]
Since
\[
\mathbb E W_\infty=\frac{\mathbb E\xi_1^2}{2\,\mathbb E\xi_1},
\]
the claim follows.
\end{proof}

\rev{The detailed quantitative and methodological comparisons are recorded in Appendix~\ref{app:benchmark_tables}.}

\begin{remark}[What the benchmark comparison shows]
\rev{The benchmark calculations show that the method is conservative at the level of constants and, in its universal form, at the level of rates. In the i.i.d.\ cases the Lorden input has the correct scale up to the universal factor $3/2$; in the mixed and dependent cases the comparison quantities remain explicit; and in the Markov-modulated case any final convergence statement must still be qualified by \eqref{cond:RD}.}
\end{remark}

\section{Conclusion}

\rev{Theorem~\ref{thm:generalized_lorden} gives a Lorden-type upper bound for the forward recurrence time in renewal-type processes with dependent, heterogeneous, and possibly mixed inter-renewal laws under Assumptions~\ref{ass:grp:min}--\ref{ass:grp:moment} and Condition~\eqref{cond:RD}. The argument is based on the generalised intensity measure and the associated comparison construction.}

\rev{Corollary~\ref{cor:generalized_lorden_coupling} records the coupling consequence only conditionally: once a stationary version, a trial-coupling scheme, and the remaining coupling parameters are available, the Lorden input yields a total-variation estimate. In this sense, the paper provides an explicit first-moment input for coupling rather than a complete unconditional convergence theorem for the whole class.}

\rev{The benchmark calculations show that, in the i.i.d.\ cases, the bound has the correct renewal scale up to the universal factor $3/2$; that mixed laws can be handled without reducing to a purely absolutely continuous setting; and that in the Markov-modulated and Pareto cases the comparison quantities remain explicit. In the Markov-modulated example the final convergence conclusion remains conditional on \eqref{cond:RD}. A natural further step would be to identify verifiable sufficient conditions for \eqref{cond:RD} and to sharpen the coupling bounds when additional structure is available.}

\section*{Acknowledgements}

\rev{The authors are grateful to Professor A.\,Yu.\ Veretennikov for valuable advice and discussions.}

\appendix

\section{\revtitle{Benchmark tables}}
\label{app:benchmark_tables}

\rev{This appendix collects the comparative benchmark tables referred to in Section~\ref{subsec:benchmark_summary}.}

Table~\ref{tab:benchmark_quantitative} summarises the quantitative calibration of the method in the benchmark models. In the i.i.d.\ benchmarks, the comparison is made both with the exact equilibrium forward-recurrence mean and, where available, with the exact total-variation convergence rate.

\begin{table}[t]
\centering
\small
\renewcommand{\arraystretch}{1.2}
\begin{tabular}{|p{2.4cm}|p{2.6cm}|p{2.8cm}|p{3.2cm}|p{3.0cm}|}
\hline
\textbf{Benchmark} &
\textbf{Exact / known Lorden input} &
\textbf{Present input \(C_{\mathrm{GL}}\)} &
\textbf{Present generic convergence bound} &
\textbf{Main message}
\\
\hline
Exponential
(Subsection~\ref{subsec:benchmark_exponential}) &
\(\mathbb E W_\infty=1/\lambda\), \quad
\(\|\Law(B_t^0)-\pi_B\|_{\mathrm{TV}}=e^{-\lambda t}\)
&
\[
C_{\mathrm{GL}}^{\mathrm{exp}}=\frac{3}{\lambda}
\]
&
Universal input alone gives
\[
\|\Law(B_t^0)-\pi_B\|_{\mathrm{TV}}
\le \frac{3}{\lambda t}.
\]
Using the exact first moment of the coupling time gives
\[
\|\Law(B_t^0)-\pi_B\|_{\mathrm{TV}}
\le \frac{1}{\lambda t},
\]
and the exact tail recovers \(e^{-\lambda t}\)
&
Correct scale \(1/\lambda\); conservativeness enters at the final coupling step, not at the renewal input itself
\\
\hline
Atom \(+\) exponential tail
(Subsection~\ref{subsec:benchmark_mixed}) &
\[
\mathbb E W_\infty
=
\frac{\mathbb E\xi^2}{2\,\mathbb E\xi}
\]
with explicit moments
&
\[
C_{\mathrm{GL}}^{\mathrm{mix}}
=
\frac{3}{2}\frac{\mathbb E\xi^2}{\mathbb E\xi}
\]
&
Explicit \(t^{-1}\)-bound once the trial parameters \(m_x,q_x\) are identified
&
The generalised intensity measure remains explicit in the presence of atoms; the correct renewal scale is preserved
\\
\hline
Markov-modulated exponential
(Subsection~\ref{subsec:benchmark_mme}) &
No simple closed exact benchmark written here; sharper matrix/spectral rates are available in principle
&
\[
C_{\mathrm{GL}}^{\mathrm{MME}}
=
\frac{2}{\phi}+\frac{\overline\lambda}{\phi^2}
\]
&
Explicit Lorden input:
\[
C_{\mathrm{GL}}^{\mathrm{MME}}
=
\frac{2}{\phi}+\frac{\overline\lambda}{\phi^2}.
\]
The full convergence-rate estimate remains conditional on \eqref{cond:RD}
&
Dependence and heterogeneity are handled explicitly; the framework and the Lorden input are fully computable, while the full convergence-rate statement remains conditional on \eqref{cond:RD}
\\
\hline
Pareto
(Subsection~\ref{subsec:benchmark_pareto}) &
\begin{align*}
\mathbb E W_\infty=\frac{b}{\alpha-2},
\\
\alpha>2
\end{align*}
&
\[
C_{\mathrm{GL}}^{\mathrm{Par}}=\frac{3b}{\alpha-2}
\]
&
\begin{align*}
\|\Law(X_t^x)-\pi\|_{\mathrm{TV}}
\le
\\
\frac{C_{\mathrm{GL}}^{\mathrm{Par}}+m_x(q_x^{-1}-1)}{t},
\end{align*}
and higher moments lead to higher-order polynomial bounds
&
Correct threshold \(\alpha>2\) for finiteness of the first Lorden input; natural polynomial regime
\\
\hline
\end{tabular}
\caption{Quantitative calibration of the present method in the benchmark models.}
\label{tab:benchmark_quantitative}
\end{table}

Table~\ref{tab:benchmark_methods} compares the present approach with standard alternatives. The point is not that the present method is uniformly sharper. Rather, the table shows where robustness is gained and where explicit constants remain available without generator or spectral analysis.

\begin{table}[t]
\centering
\small
\renewcommand{\arraystretch}{1.2}
\begin{tabular}{|p{2.5cm}|p{3.2cm}|p{3.2cm}|p{4.0cm}|}
\hline
\textbf{Benchmark} &
\textbf{Generator / spectral route} &
\textbf{Lyapunov route} &
\textbf{Present method}
\\
\hline
Exponential &
Exact and explicit; recovers the sharp exponential rate \(e^{-\lambda t}\)
&
Available but unnecessary in this benchmark
&
Generic first-moment route gives a conservative polynomial bound; refined coupling recovers the exact exponential rate
\\
\hline
Atom \(+\) exponential tail &
Possible in principle, but explicit generator-based estimates become less natural in the presence of atoms and typically require model-specific analysis
&
Possible in principle, but requires a model-specific drift construction adapted to the mixed law
&
The generalised intensity measure is explicit; the Lorden input is obtained in closed form directly from the mixed law
\\
\hline
Markov-modulated exponential &
Available and potentially sharper via matrix-analytic or spectral methods, but depends on the full transition structure of the modulating chain
&
Available after Markov augmentation, but explicit constants are model-specific
&
Fully explicit and generator-free at the level of the comparison framework and Lorden input
\\
\hline
Pareto &
A generator exists for the age process in the i.i.d.\ benchmark, but quantitative bounds are naturally polynomial and model-specific
&
Polynomial Lyapunov drifts are conceivable, but explicit constants are not automatic
&
The comparison intensity is explicit, the threshold \(\alpha>2\) is identified correctly, and the method naturally produces polynomial regimes
\\
\hline
\end{tabular}
\caption{Methodological comparison across the benchmark models.}
\label{tab:benchmark_methods}
\end{table}









\bibliographystyle{unsrt}
\bibliography{sn-bibliography}

\end{document}